%
%
%
%
\hsize=5in
\baselineskip=12pt
\vsize=20.4cm
\parindent=.5cm
\predisplaypenalty=0
\hfuzz=2pt
\frenchspacing
\def\latexfmt{latex}
\ifx\fmtname\latexfmt\else
%
%
\input amssym.def
\def\titlefonts{\baselineskip=1.44\baselineskip
	\font\titlef=cmbx12
	\titlef
	}
\font\ninerm=cmr9
\font\ninebf=cmbx9
\font\ninei=cmmi9
\skewchar\ninei='177
\font\ninesy=cmsy9
\skewchar\ninesy='60
\font\nineit=cmti9
\def\reffonts{\baselineskip=0.9\baselineskip
	\textfont0=\ninerm
	\def\rm{\fam0\ninerm}%
	\textfont1=\ninei
	\textfont2=\ninesy
	\textfont\bffam=\ninebf
	\def\bf{\fam\bffam\ninebf}%
	\def\it{\nineit}%
	}
%
%
\def\frontmatter{\vbox{}\vskip3cm\bgroup
	\leftskip=0pt plus1fil\rightskip=0pt plus1fil
	\parindent=0pt
	\parfillskip=0pt
	\pretolerance=10000
	}
\def\endfrontmatter{\egroup\bigskip}
\def\title#1{{\titlefonts#1\par}}
\def\author#1{\bigskip#1\par}
\def\address#1{\bigskip{\reffonts\it#1}}
\def\email#1{\bigskip{\reffonts{\it E-mail: }\rm#1}}
\def\thanks#1{\footnote{}{\reffonts\rm\noindent#1\hfil}}
\def\section#1\par{\ifdim\lastskip<\bigskipamount\removelastskip\fi
	\penalty-250\bigskip
	\vbox{\leftskip=0pt plus1fil\rightskip=0pt plus1fil
	\parindent=0pt
	\parfillskip=0pt
	\pretolerance=10000{\bf#1}}\nobreak\medskip
	}

\fi
\def\proclaim#1. {\medbreak\bgroup{\noindent\bf#1.}\ \it}
\def\endproclaim{\egroup
	\ifdim\lastskip<\medskipamount\removelastskip\medskip\fi}
\def\item#1 #2\par{\ifdim\lastskip<\smallskipamount\removelastskip\smallskip\fi
	{\rm#1}\ #2\par\smallskip}
\def\Proof#1. {\ifdim\lastskip<\medskipamount\removelastskip\medskip\fi
	{\noindent\it Proof.}\ }
\def\endproof{\hfill\hbox{$\square$}\medskip}
\def\Remark. {\ifdim\lastskip<\medskipamount\removelastskip\medskip\fi
	{\noindent\bf Remark.}\quad}
\def\Remarks. {\ifdim\lastskip<\medskipamount\removelastskip\medskip\fi
	{\noindent\bf Remarks.}\quad}
\def\endremark{\medskip}
%
%
\newcount\citation
\newtoks\citetoks
\def\citedef#1\endcitedef{\citetoks={#1\endcitedef}}
\def\endcitedef#1\endcitedef{}
\def\citenum#1{\citation=0\def\curcite{#1}%
	\expandafter\checkendcite\the\citetoks}
\def\checkendcite#1{\ifx\endcitedef#1?\else
	\expandafter\lookcite\expandafter#1\fi}
\def\lookcite#1 {\advance\citation by1\def\auxcite{#1}%
	\ifx\auxcite\curcite\the\citation\expandafter\endcitedef\else
	\expandafter\checkendcite\fi}
\def\cite#1{\makecite#1,\cite}
\def\makecite#1,#2{[\citenum{#1}\ifx\cite#2]\else\expandafter\clearcite\expandafter#2\fi}
\def\clearcite#1,\cite{, #1]}
%
%
\def\references{\section References\par
	\bgroup
	\parindent=0pt
	\reffonts
	\rm
	\frenchspacing
	\setbox0\hbox{99. }\leftskip=\wd0
	}
\def\endreferences{\egroup}
\newtoks\nextauth
\newif\iffirstauth
\def\checkendauth#1{\ifx\endauth#1%
		\iffirstauth\the\nextauth
		\else{} and \the\nextauth\fi,
	\else\iffirstauth\the\nextauth\firstauthfalse
		\else, \the\nextauth\fi
		\expandafter\auth\expandafter#1\fi
	}
\def\auth#1,#2;{\nextauth={#1 #2}\checkendauth}
\newif\ifinbook
\newif\ifbookref
\def\nextref#1 {\par\hskip-\leftskip
	\hbox to\leftskip{\hfil\citenum{#1}.\ }%
	\initnextref}
\def\initnextref{\bookreffalse\inbookfalse\firstauthtrue\ignorespaces}
\def\paper#1{{\it#1},}
\def\book#1{\bookreftrue{\it#1},}
\def\journal#1{#1}
\def\bookseries#1{#1,}
\def\Vol#1{\ifbookref Vol. #1,\else{\bf#1}\fi}
\def\publisher#1{#1,}
\def\Year#1{\ifbookref#1\else(#1)\fi}
\def\Pages#1{#1.}
%
%
\newsymbol\square 1003
\let\Rar\Rightarrow
\let\Lrar\Leftrightarrow
\let\ot\otimes
\let\sbs\subset
\newsymbol\rightrightarrows 1313
\newsymbol\smallsetminus 2272
\let\setm\smallsetminus
\def\chr{\mathop{\fam0 char}\nolimits}
\def\Com{\mathop{\fam0 Com}\nolimits}
\def\Der{\mathop{\fam0 Der}\nolimits}
\def\Hom{\mathop{\fam0 Hom}\nolimits}
\def\id{\mathop{\fam0 id}\nolimits}
\def\im{\mathop{\fam0 im}}
\def\Im{\mathop{\fam0 Im}}
\def\rk{\mathop{\fam0 rk}\nolimits}
\def\Spec{\mathop{\fam0 Spec}}
\def\St{\mathop{\fam0 St}}
\def\tr{\mathop{\fam0 tr}\nolimits}
\def\textsum{\textstyle\sum\limits}
\def\mapr#1{{}\mathrel{\smash{\mathop{\longrightarrow}\limits^{#1}}}{}}
\def\lmapr#1#2{{}\mathrel{\smash{\mathop{\count0=#1
  \loop
    \ifnum\count0>0
    \advance\count0 by-1\smash{\mathord-}\mkern-4mu
  \repeat
  \mathord\rightarrow}\limits^{#2}}}{}}
\def\mapd#1#2{\llap{$\vcenter{\hbox{$\scriptstyle{#1}$}}$}\big\downarrow
  \rlap{$\vcenter{\hbox{$\scriptstyle{#2}$}}$}}
\def\diagram#1{\vbox{\halign{&\hfil$##$\hfil\cr #1}}}
\def\diagramskip{\noalign{\smallskip}}
\let\al\alpha
\let\be\beta
\let\ga\gamma
\let\de\delta
\let\ep\varepsilon
\let\io\iota
\let\la\lambda
\let\ph\varphi
\let\si\sigma
\let\De\Delta
\newsymbol\diagdown 201F
\def\tsear#1#2{\setbox0\hbox{$\diagdown$}\setbox1\hbox{$\searrow$}%
	\setbox2\hbox{$\scriptstyle#2$}%
	\setbox3\hbox{\kern1pt\raise8.5pt\copy0\kern-\wd2\raise-3pt\box2\kern-0.8pt\box0\kern-1.7pt\raise-8pt\box1}%
	\ht3=0pt\dp3=0pt
	{}\mathrel{\raise#1\box3}{}}
\def\a{\overline a}
\def\A{\overline A}
\def\g{\overline\ga}
\def\m{{\frak m}}
\def\M{{\cal M}}
\def\n{{\frak n}}
\def\p{{\frak p}}
\def\q{{\frak q}}
\def\Z{{\Bbb Z}}

\citedef
Ar86
Ar87
BouA
Bou
Br
Chase
De
Doi85
Doi89
Doi90
Ea
Hab75
Hoch54
Ho71
Kop93
Kr81
Lar95
Mas92
Ma
Mo
Mum65
Mum70
Na
Par71
Sch90
Sk02
Sw
Ul82
Wat
Zhu96
\endcitedef
%
%
\frontmatter

\title{Invariants of finite Hopf algebras}
\author{Serge Skryabin}
\address{Chebotarev Research Institute,
Kazan, Russia}
\address{Current address: Mathematisches Seminar, University of Hamburg,\break
Bundesstr. 55, 20146 Hamburg, Germany}
\email{fm1a009@math.uni-hamburg.de}
\thanks{The author gratefully acknowledges the support of the Alexander
von Humboldt Foundation through the program of long term cooperation}

\endfrontmatter

\section
Introduction

This paper extends classical results in the invariant theory of finite groups
and finite group schemes to the actions of finite Hopf algebras on commutative
rings. Suppose that $H$ is a finite dimensional Hopf algebra and $A$ a
commutative algebra, say over a field $K$. Let $\de:A\to A\ot H$ be an algebra
homomorphism which makes $A$ into a right $H$-comodule. In this case $A$ is
called an $H$-comodule algebra. The coaction of $H$ on $A$ corresponds to an
action of the dual Hopf algebra $H^*$. For technical reasons all results in
this paper are formulated in terms of coactions. The situation where $H$ is
commutative can be described geometrically by giving an action of the finite group
scheme $G=\Spec H$ on the scheme $X=\Spec A$. The subring of $G$-invariants
$A^G\sbs A$ represents then the quotient scheme $X/G$. A fact of fundamental
importance states that $A$ is an integral extension of $A^G$. In case of ordinary
finite groups this classical result goes back to the work of E.~Noether. The
question of whether a similar assertion is true for a noncommutative $H$ was posed
by Montgomery \cite{Mo, 4.2.6}. Shortly afterwards Zhu \cite{Zhu96} succeeded in
verifying two special cases and constructing a counterexample in general.

The first objective of the present article is to investigate the integrality
of $A$ over $A^H$ more thoroughly (here $A^H$ stands for the invariants of the
given coaction). A part of the argument given in \cite{De} and \cite{Mum70} for
commutative $H$ carries over without problems. Since $A$ is commutative, $A\ot H$
can be regarded as an $A$-algebra with respect to the action on the first tensorand.
Since $A\ot H$ is free of finite rank over $A$, with each element $u\in A\ot H$ one
can associate its characteristic polynomial $P_{A\ot H/A}(u,t)\in A[t]$. Letting
$$
P_{A\ot H/A}(\de a,t)=\textsum_{i=0}^nc_it^i
$$
for $a\in A$, where $n=\dim H$, one has
$c_n=1$ and $\sum_{i=0}^nc_ia^i=0$. If $H$ is commutative then
$c_0,\ldots,c_n\in A^H$, and the desired integrality is achieved. At this
point the noncommutativity of $H$ spoils the game altogether. It is
nevertheless shown in Theorem 2.5 that $c_0,\ldots,c_n\in A^H$ for any $H$
provided that $A$ is $H$-reduced, that is, $A$ has no nonzero $H$-costable nil
ideals. Crucial for the proof is the freeness of finite dimensional Hopf
algebras over commutative right coideal subalgebras. Because of nice functorial
properties of characteristic polynomials, it is possible to pass to finite
dimensional Hopf algebras $E\ot H$ over various fields $E$ where the
characteristic polynomials can be computed using the above mentioned freeness.

The integrality of $A$ over $A^H$ is proved in Proposition 2.7 and Theorem 6.2
not only for an $H$-reduced $A$ but also in other situations. It always holds
when either $\chr K\ne0$, or $H$ is cosemisimple, or there exists a total integral
$H\to A$ in the sense of Doi \cite{Doi85}. In Zhu's paper \cite{Zhu96, Th.~2.1, Cor.~3.3}
the integrality was proved under the assumptions that either $\chr K$ does not divide
$\dim H$ and $H$ is involutory or $\chr K\ne0$ and $H^*$ has a
cocommutative coradical. In the first of these two results $H$ is necessarily semisimple
and cosemisimple. At the same time it is still not known if every semisimple
Hopf algebra is involutory when $\chr K\ne0$ \cite{Lar95}.

The invariance of characteristic polynomials enables one to investigate the
map $\Spec A\to\Spec A^H$ between the prime spectra of rings $A$ and $A^H$
much in the same spirit as was done in case of commutative $H$. It is shown in Theorem
3.3 that this map has finite fibers, is open and satisfies the going-down
property.

With each $\p\in\Spec A$ we associate the orbital subalgebra $O(\p)\sbs k(\p)\ot H$
and the stabilizer subalgebra $\St(\p)\sbs k(\p)\ot H^*$ where $k(\p)$ denotes the
residue field of the local ring $A_\p$. A prime $\p$ will be called
$H$-regular if $\dim_{k(\p)}O(\p)$ does not change in a suitable neighbourhood
of $\p$ in $\Spec A$. The assumption that the algebra $A$ is $H$-reduced and
all prime ideals of $A$ are $H$-regular implies that $A$ is a finitely
generated projective $A^H$-module and the $H$-costable ideals of $A$
correspond bijectively to the ideals of $A^H$. This is the content of Theorem
4.3 which extends my earlier work \cite{Sk02}. The conclusions of this theorem
were known to be true when $A$ is an $H$-Galois extension of $A^H$.

The purpose of section 5 is to show the existence of a total integral $H\to A$
for an $H$-comodule algebra with semisimple stabilizer subalgebras $\St(\p)$.
This is one of cases in which the canonical maps $A^H\to(A/I)^H$ are
surjective for all $H$-costable ideals $I$ of $A$. In general, if the latter
property is fulfilled, $A$ will be called weakly reductive with respect to
coaction of $H$. This property is essential in the final major result of this paper.
Theorem 6.2 states that $A^H$ is Cohen-Macaulay whenever $A$ is Cohen-Macaulay
and weakly reductive. This generalizes a result of Hochster and Eagon \cite{Ho71,
Prop.~13} according to which the Cohen-Macaulay property of a commutative ring $A$
descends to the subring of invariants of a finite group $G$ acting on $A$
by automorphisms provided that the order of $G$ is invertible in $A$. A point
of interest here is that $A$ is always weakly reductive when $\chr K=0$, even
if $H$ is not cosemisimple. Theorem 6.2 also gives sufficient conditions for
$A$ to be a finitely generated $A^H$-module.

Summarizing, the invariants of a noncommutative $H$ behave as decently as in
the commutative case provided that one is going to be content with the
$H$-reduced algebras $A$. The assumption that $K$ is a field was made in this
introduction only to illustrate the results in the simplest case. Working over
a commutative base ring involves no significant complications.

\section
1. Preliminaries

We fix a few notations for the whole paper. Assume that $K$ is a commutative
ring, $H$ is a Hopf algebra over $K$ whose underlying $K$-module is
finitely generated and projective, and $A$ is a commutative right $H$-comodule algebra.
We will consider mainly right comodules, and the prefix ``right" will
be omitted. All tensor products, when the base ring is not indicated, are
taken over $K$. The comultiplication, counit and antipode in $H$ are
denoted as $\De$, $\ep$ and $\si$, respectively. We write symbolically
$\De h=\sum_{(h)}h'\ot h''$ and $(\De\ot\id)\De h=\sum_{(h)}h'\ot h''\ot h'''$
for $h\in H$.

For every $H$-comodule $M$ denote by $\de:M\to M\ot H$ the corresponding
structure map. In particular, $\de:A\to A\ot H$ is a homomorphism of unital
algebras. Denote by $M^H\sbs M$ the $K$-submodule of invariants. Thus
$M^H=\{v\in M\mid\de v=v\ot1\}$. Clearly $A^H$ is a subalgebra of $A$. Next,
$\M_A^H$ stands for the category of right $(H,A)$-Hopf modules \cite{Doi85}.
The objects of $\M_A^H$ are right $A$-modules equipped with a right
$H$-comodule structure such that
$$
\de(va)=\de(v)\de(a)\qquad\qquad{\fam0for\ all\ }v\in M
{\fam0\ and\ }a\in A.\eqno(1.1)
$$
Here and later we regard $M\ot H$ as a right $A\ot H$-module by means of the
operation $(v\ot g)(a\ot h)=va\ot gh$ for $v\in M$, $a\in A$ and $g,h\in H$.
The morphisms in $\M_A^H$ are maps which are compatible with both
structures.

\proclaim
Lemma 1.1.
Suppose that $S\sbs A^H$ is any multiplicatively closed subset. Then
$S^{-1}A$ is an $H$-comodule algebra and $S^{-1}M\in\M_{S^{-1}A}^H$ for every
$M\in\M_A^H$. Moreover, one has $(S^{-1}M)^H\cong S^{-1}(M^H)$.
\endproclaim

\Proof.
Let us regard $M\ot H$ as an $A^H$-module by means of the action on the first
tensorand. The two maps $\de,\io:M\to M\ot H$ where $\io(v)=v\ot1$ for $v\in
M$ are homomorphisms of $A^H$-modules. Hence $\de$ extends to a map
$S^{-1}M\to S^{-1}M\ot H$ which makes $S^{-1}M$ into an $H$-comodule. One
checks easily that (1.1) is fulfilled for the action of $S^{-1}A$ on
$S^{-1}M$. The assertion about $(S^{-1}M)^H$ is obtained by applying the
localization functor $S^{-1}(?)$ to the exact sequence of $A^H$-modules
$$
0\mapr{}M^H\mapr{}M\lmapr3{\de-\io}M\ot H.
\vadjust{\vskip-16pt}
$$
\endproof

This lemma will be used in two cases. If $S=\{s^i\mid i\ge0\}$ where $s\in
A^H$ then the localizations are denoted as $A_s$ and $M_s$. If $S=A^H\setm\q$
where $\q\in\Spec A^H$ then the notations are $A_\q$ and $M_\q$.

An ideal $I$ of $A$ is called {\it$H$-costable} if $\de(I)\sbs I\ot H$. This
makes sense since $I\ot H$ is embedded into $A\ot H$ by the $K$-projectivity
of $H$. If $I$ is $H$-costable then $\de:A\to A\ot H$ induces a ring
homomorphism $A/I\to A/I\ot H$ which makes $A/I$ into an $H$-comodule algebra.

If $\q\in\Spec K$ then $H_\q=K_\q\ot H$ is a free module over the local ring
$K_\q$. Denote by $\rk_\q H$ its rank. If $\rk_\q H$ does not depend on $\q$,
then $H$ is $K$-projective of {\it constant rank}.

\proclaim
Lemma 1.2.
One has $K\cong\prod K_i$, $\,H\cong\prod H_i$, $\,A\cong\prod A_i$ where the
products are taken over a finite set of indices, each $K_i$ is a commutative
ring, $H_i$ is a Hopf algebra over $K_i$ which is $K_i$-projective of constant
rank, and $A_i$ is an $H_i$-comodule algebra. Moreover, $A^H\cong\prod A_i^{H_i}$.
\endproclaim

\Proof.
For each $i>0$ the subset $X_i=\{\q\in\Spec K\mid\rk_\q H=i\}$ is open in
$\Spec K$ \cite{Bou, Ch.~II, \S5, Th.~1}. Now $\Spec K$ is a finite disjoint
union of these subsets $X_i$, where we remove the empty subsets. By
\cite{Bou, Ch.~II, \S4, Prop.~15} $K$ contains a family of orthogonal
idempotents $\{e_i\}$ such that $\sum e_i=1$ and $X_i=\{\q\in\Spec
K\mid1-e_i\in\q\}$ for each $i$. If $V$ is any $K$-module then $V\cong\prod
V_i$ where $V_i=V/(1-e_i)V$. This gives the first three isomorphisms of the
lemma. Similarly, $A\ot H\cong\prod(A_i\ot H_i)$ and the map $\de:A\to A\ot
H$ is compatible with the cartesian products. Hence the assertion about $A^H$.
\endproof

As is well known, the dual $H^*=\Hom_K(H,K)$ of $H$ is a Hopf algebra in a
natural way. The $H$-comodule structures are in a bijective correspondence
with the left $H^*$-module structures \cite{Par71, Prop.~1}. The operator
giving the action of $\xi\in H^*$ on an $H$-comodule $M$ can be written as the
composite
$$
L_\xi:M\mapr\de M\ot H\lmapr4{\id\ot\xi}M\ot K\cong M.\eqno(1.2)
$$
One sees that $M^H=\{v\in M\mid(H^*)^+v=0\}$ where $(H^*)^+=\{\xi\in
H^*\mid\xi(1)=0\}$. We may regard $H$ as a right $H$-comodule via $\De$.

\proclaim
Lemma 1.3.
$H$ is a finitely generated projective $H^*$-module with respect to the
corresponding action of $H^*$.
\endproclaim

\Proof.
By \cite{Par71, Prop.~2, Lemma 2, Prop.~3} with $H$ and $H^*$
interchanged the set of left integrals $P(H)=\{x\in H\mid hx=\ep(h)x$ for all
$h\in H\}$ is a rank one projective $K$-module and the map $H^*\ot P(H)\to H$
defined by the rule $\xi\ot x\mapsto L_\xi(x)$ is bijective. Clearly $H^*\ot
P(H)$ is a finitely projective $H^*$-module.
\endproof

A reduction of the base ring to a field is one of our tools. We will need
several facts from earlier work. Assume until the end of this section
that {\it K is a field} and {\it H is a finite dimensional Hopf algebra}. A
{\it left} (respectively, {\it right}) {\it coideal subalgebra} $B\sbs H$ is a
subalgebra satisfying $\De(B)\sbs H\ot B$ (respectively, $\De(B)\sbs B\ot H$).
If $H$ is cocommutative then the two inclusions are equivalent to one another,
and so every left coideal subalgebra of $H$ is a Hopf subalgebra. The next
result is due to Masuoka \cite{Mas92, (2.1)}:

\proclaim
Proposition 1.4.
Suppose that $B\sbs H$ is a right coideal subalgebra. Then $B$ is Frobenius
if and only if $H$ is left and right $B$-free.
\endproclaim

\proclaim
Proposition 1.5.
Every commutative right coideal subalgebra $B\sbs H$ is Frobenius, and so
$H$ is a free $B$-module from either side.
\endproclaim

This is contained in Koppinen's paper \cite{Kop93, Cor.~2.5} where the result
is attributed to the referee.

\proclaim
Proposition 1.6.
Suppose that $B\sbs H$ is a Frobenius right coideal subalgebra and
$B^+=\{b\in B\mid\ep(b)=0\}$. Then{\rm:}

\item(i)
$C=(H/HB^+)^*$ is a Frobenius left coideal subalgebra of $H^*${\rm;}

\item(ii)
$B\cong(H^*/C^+H^*)^*$ where $C^+=\{\xi\in C\mid\xi(1)=0\}${\rm;}

\item(iii)
The category $\M_B^H$ is equivalent to the category of $C$-modules.

\endproclaim

This is a version of Masuoka's result \cite{Mas92, (2.10)} where a similar
relationship between the right coideal subalgebras in both $H$ and $H^*$ is
described. One connects the two formulations by taking the opposite
multiplications in $A$ and $H$ and applying \cite{Mas92, (1.1), (1.2)}.

\section
2. Characteristic polynomials and integrality over invariants

Suppose that $A$ is any commutative ring and $U$ is an associative
$A$-algebra such that $U$ is projective of constant rank $n$ as an
$A$-module. In this case one can define the {\it norm} $N_{U/A}(u)\in A$ for
every $u\in U$ \cite{BouA, Ch.~III, \S9}, \cite{Bou, Ch.~II, \S5, Exercise 9}. In fact, if
$U$ has a basis $e_1,\ldots,e_n$ over $A$, then $ue_j=\sum_{i=1}^na_{ij}e_i$
with coefficients in $A$, and $N_{U/A}(u)$ is the determinant of the square
matrix with entries $a_{ij}$, $\,1\le i,j\le n$. In general one can pass to
localizations $U_s$ and $A_s$ with respect to multiplicatively closed
subsets $\{s^i\mid i\ge0\}$ where $s\in A$. There exist a finite number of
elements $s_1,\ldots,s_m\in A$ which generate the whole $A$ as an ideal and
such that $U_{s_i}$ is free of rank $n$ over $A_{s_i}$ for each
$i=1,\ldots,m$. If $u_i$ denotes the image of $u$ in $U_{s_i}$, then
$N_{U_{s_i}/A_{s_i}}(u_i)\in A_{s_i}$ is defined, and any two of these norms
corresponding to a pair of indices $i,j$ have the same image in $A_{s_is_j}$.
Then this collection of elements in localizations can be glued to an element
$N_{U/A}(u)\in A$ using the well known exact sequence
$$
0\to A\to\prod_iA_{s_i}\rightrightarrows\prod_{i,j}A_{s_is_j}.
$$
The {\it characteristic polynomial} of $u$ is obtained by adjoining an indeterminate
$t$ as
$$
P_{U/A}(u,t)=N_{U[t]/A[t]}(t-u)\in A[t].
$$
In particular, $(-1)^nN_{U/A}(u)$ is the coefficient of $t^0$ in
$P_{U/A}(u,t)$. The well known properties of the characteristic polynomials and
the norm are listed below:

\item(P1)
$P_{U/A}(u,u)=0$ (the substitution of $u$ for $t$);

\item(P2)
$P_{V/A}(\ph u,t)=P_{U/A}(u,t)$ when $\ph:U\to V$ is an isomorphism of
$A$-algebras;

\item(P3)
$P_{B\ot_AU/B}(1\ot u,t)=\ga^tP_{U/A}(u,t)$ where $\ga:A\to B$ is a
homomorphism of commutative rings and $\ga^t:A[t]\to B[t]$ is its extension
such that $t\mapsto t$;

\item(P4)
If $V\subset U$ is an $A$-subalgebra such that $V$ is projective as an
$A$-module and $U$ is free, say of rank $r$, as a left $V$-module, then
$P_{U/A}(u,t)=P_{V/A}(u,t)^r$ for every $u\in V$;

\item(P5)
If $a\in A$ then $P_{U/A}(a,t)=(t-a)^n$;

\item(P6)
If $A$ is a field and $\dim_AU=n$, then the equality $P_{U/A}(u,t)=t^n$ is a
necessary and sufficient condition for $u\in U$ to be nilpotent;

\item(P7)
An element $u$ is invertible in $U$ if and only if $N_{U/A}(u)$ is invertible
in $A$.

Suppose further that $A$ is a commutative $H$-comodule algebra as was
specified in section 1.  Given a ring homomorphism $\al:A\to E$ into a field
$E$, define
$$
A_\al=(E\ot1)\cdot\de_\al(A)\subset E\ot H\eqno(2.1)
$$
where $\de_\al$ is the composite
$$
A\mapr\de A\ot H\lmapr4{\al\ot\id}E\ot H.\eqno(2.2)
$$

\proclaim
Lemma 2.1.
$A_\al$ is a right coideal subalgebra of the Hopf algebra $E\ot H$ over $E$.
If $\De_E$ denotes the map $\id\ot\De:E\ot H\to E\ot H\ot H$, then
$$
(\de_\al\ot\id)\circ\de=\De_E\circ\de_\al\,.\eqno(2.3)
$$
\endproclaim

\Proof.
Since $\de$ is a ring homomorphism, so too is $\de_\al$, whence $A_\al$ is a
subalgebra of $E\ot H$. Formula (2.3) is seen from the commutative diagram
$$
\diagram{
A&\hidewidth\lmapr7\de\hidewidth&A\ot H&\hidewidth\lmapr{10}{\al\ot\id}\hidewidth&E\ot H\cr
\diagramskip
\mapd{}\de&&\mapd{}{\id\ot\De}&&\mapd{}{\id\ot\De}\cr
\diagramskip
A\ot H&\lmapr3{\de\ot\id}&A\ot H\ot H&\lmapr6{\al\ot\id\ot\id}&E\ot H\ot H.\cr
}
$$
It follows from (2.3) that $\De_E(A_\al)\subset A_\al\ot H$. Under the
canonical identification $(E\ot H)\ot_E(E\ot H)\cong E\ot H\ot H$, the
comultiplication in $E\ot H$ is precisely $\De_E$. Hence $A_\al$
is a right coideal.
\endproof

\Remark.
One can also view $E\ot H$ as an $H$-comodule algebra with respect to $\De_E$
and $A_\al$ as its $H$-comodule subalgebra. Formula (2.3) says that
$\de_\al:A\to E\ot H$ is a homomorphism of $H$-comodule algebras.
\endremark

\proclaim
Lemma 2.2.
Suppose that $B\subset C$ is an extension of commutative rings such that $B$
is Artinian and $C$ is $B$-projective. If $M$ is a $C$-module such that $M$
is $B$-projective and $M/\n M$ is free of rank $r$ over $C/\n C$ for every
maximal ideal $\n$ of $B$ where $r$ does not depend on $\n$, then $M$ is
free of rank $r$ over $C$.
\endproclaim

\Proof.
Let $J$ be the Jacobson radical of $B$. Then $B/J\cong\prod B/\n$, the product over
the finitely many maximal ideals of $B$. Furthermore, $C/JC\cong\prod C/\n C$ and
$M/JM\cong\prod M/\n M$ \cite{Bou, Ch.~II, \S1, Prop.~6}. It follows from the hypotheses
of the lemma that $M/JM$ is a free $C/JC$-module of rank $r$. Let $F$ be a free $C$-module of rank $r$.
There exists then a homomorphism of $C$-modules $\ph:F\to M$ which induces an isomorphism
$F/JF\cong M/JM$. Since $J$ is a nilpotent ideal of $B$, we conclude that $\ph$ is surjective
\cite{Bou, Ch.~II, \S3, Cor.~1 to Prop.~4}. Put $N=\ker\ph$. By the $B$-projectivity of $M$
the exactness of sequence $0\to N\to F\to M\to0$ is preserved when passing to reductions
modulo $J$. Hence $N/JN=0$, and $N=0$ by Nakayama's Lemma. Thus $\ph$ is an isomorphism.
\endproof

The ring $A\ot H$ contains $A\ot1$ in its center. Hence we may view $A\ot H$
as an $A$-algebra using the canonical homomorphism $A\to A\ot1$. Since $H$
is finitely projective $K$-module, so too is $A\ot H$ as an
$A$-module. If $H$ is a projective $K$-module of rank $n$, the $A$-module
$A\ot H$ has rank $n$ as well. In this case the characteristic polynomial
$P_{A\ot H/A}(u,t)\in A[t]$ is defined for every $u\in A\ot H$.

\proclaim
Proposition 2.3.
Suppose that $K$ is a field. If $b\in B$ where $B$ is a commutative right
coideal subalgebra of $H$, then $P_{B\ot H/B}(\De b,t)=P_{H/K}(b,t)$.
\endproclaim

The conclusion of this proposition means, in particular, that the first
polynomial has all coefficients in $K\subset B$.

\Proof.
We regard $B\ot H$ as a $B$-algebra in accordance with the previous
explanation. By Proposition 1.5 $H$ is free over $B$ with respect to
the action by left multiplications. Let $r=\rk_BH$. Put
$C=(B\ot1)\cdot\De(B)\subset B\ot H$. Since $\De:B\to B\ot H$ is a ring
homomorphism, $C$ is a commutative $B$-subalgebra of $B\ot H$. We will prove
first that $B\ot H$ is a free $C$-module of rank $r$ with respect to the
action by left multiplications.

Define $\Phi:B\ot B\to B\ot H$ by $\Phi(a\ot b)=(a\ot1)\cdot\De b$ for
$a,b\in B$. Then $\Phi$ is a ring homomorphism and $\Im\Phi=C$. Since $B$ is
commutative and $H$ is a free $B$-module of finite rank, $B$ is a
$B$-module direct summand of $H$ \cite{Bou, Ch.~II, \S5, Excersize 4}. Let $\psi:H\to B$
be a retraction, so that $\psi(1)=1$ and $\psi(bh)=b\psi(h)$ for all $b\in B$ and $h\in H$.
Define $\Psi:B\ot H\to B\ot B$ by $\Psi(a\ot h)=\sum_{(h)}a\psi(\si h')\ot h''$. Then
$$
\Psi(\De b)=\textsum_{(b)}\,b'\psi(\si b'')\ot b'''
=\textsum_{(b)}\,\psi(b'\cdot\si b'')\ot b'''
=\textsum_{(b)}\,\ep(b')\psi(1)\ot b''=1\ot b.
$$
for $b\in B$. This shows that $\Psi\circ\Phi=\id$. In particular, $\Phi$ is a ring
isomorphism of $B\ot B$ onto $C$. We have also $B\ot H=C\oplus\ker\Psi$.
Both $\Phi$ and $\Psi$ are $B$-linear with respect to the actions of
$B$ on the first tensorands. Then $C$ is a $B$-module direct summand of
$B\ot H$. Let $\n$ be a maximal ideal of $B$ and $E=B/\n$ its residue field.
The embedding of $C$ into $B\ot H$ induces an injective homomorphism of
$E$-algebras
$$
C/\n C\to(B\ot H)/\n(B\ot H)\cong E\ot H.
$$
The image of the latter coincides with $B_\be$ where $\be:B\to E$ is the
canonical projection and $B_\be$ is defined as in (2.1). By Lemma 2.1
$B_\be$ is a right coideal subalgebra of $E\ot H$, and by Proposition 1.5
$E\ot H$ is free over $B_\be$. Now the reduction of $\Phi$ modulo
$\n$ gives an $E$-linear bijection $E\ot B\cong C/\n C\cong B_\be$. Hence
$\dim_EB_\be=\dim_KB$, and so
$$
\rk_{B_\be}E\ot H=\dim_E(E\ot H)/\dim_EB_\be=\dim_KH/\dim_KB=r.
$$
We see that the hypotheses of Lemma 2.2 are fulfilled for $M=B\ot H$,
yielding the desired freeness.

Let $b\in B$, and so $\De b\in C$. Applying successively (P4), (P2), (P3)
and (P4) again, we find
$$
P_{B\ot H/B}(\De b,t)=P_{C/B}(\De b,t)^r=P_{B\ot B/B\ot1}(1\ot b,t)^r
=P_{B/K}(b,t)^r=P_{H/K}(b,t).
\vadjust{\vskip-8pt}
$$
\endproof

\proclaim
Lemma 2.4.
Suppose that $\al:A\to E$ is a ring homomorphism into a field $E$. If
$b=\de_\al(a)$ where $a\in A$, then
$$
P_{A_\al\ot H/A_\al}(\De_Eb,t)=P_{E\ot H/E}(b,t).\eqno(2.4)
$$
\endproclaim

\Proof.
Apply Proposition 2.3 to the Hopf algebra $E\ot H$ over the field $E$ and
its right coideal subalgebra $B=A_\al$ (see Lemma 2.1).
\endproof

When $H$ is $K$-projective of constant rank, we say that $A$ has {\it
invariant characteristic polynomials} if $P_{A\ot H/A}(\de a,t)$ has all
coefficients in $A^H$ for every $a\in A$. When the rank of $H$ is not
constant, we say that $A$ has invariant characteristic polynomials if so do
the $H_i$-comodule algebras $A_i$ in the cartesian product decomposition of
Lemma 1.2. In case of a commutative $H$ the invariance of characteristic
polynomials for any $A$ was proved in \cite{De, Ch.~III, \S12} and
\cite{Mum70, Ch.~III, \S2}. If $H$ is not commutative, additional assumptions
are needed. The largest $H$-costable nil ideal of $A$ will be called the
{\it$H$-radical} of $A$. We say that $A$ is {\it$H$-reduced} if its
$H$-radical is equal to $(0)$.

\proclaim
Theorem 2.5.
If $A$ is $H$-reduced or, more generally, if there exists a homomorphism of
commutative $H$-comodule algebras $\ph:A'\to A$ such that $A'$ is
$H$-reduced and $A=\ph(A')A^H$ then $A$ has invariant characteristic
polynomials.
\endproclaim

\Proof.
Consider first the case where $A$ is $H$-reduced. The $H_i$-comodule algebras
$A_i$ in Lemma 1.2 are then $H_i$-reduced too. So we may assume without loss
of generality that $H$ is $K$-projective of constant rank.

Suppose that $\al:A\to E$ is a ring homomorphism into a field $E$. Define
$\de_\al$ as in (2.2), and let $\la_\al=\io\circ\al$ where $\io:E\to E\ot H$
is given by $c\mapsto c\ot1$. Both $\de_\al$ and $\la_\al$ take values in
$A_\al$. We first prove that
$$
\de_\al^tP_{A\ot H/A}(\de a,t)=\la_\al^tP_{A\ot H/A}(\de a,t)
\quad{\fam0in\ }A_\al[t].\eqno(2.5)
$$
If $\ga:A\to A_\al$ is any ring homomorphism then
$$
\ga^tP_{A\ot H/A}(\de a,t)
=P_{A_\al\ot H/A_\al}\bigl((\ga\ot\id)\de a,t\bigr).
$$
by (P3) since $A_\al\ot_A(A\ot H)\cong A_\al\ot H$ and the element $1\ot u$,
where $u\in A\ot H$, corresponds under this isomorphism to $(\ga\ot\id)u\in A_\al\ot H$.
In case $\ga=\de_\al$, denoting $b=\de_\al(a)$ and using (2.3), (2.4), we obtain
$$
\de_\al^tP_{A\ot H/A}(\de a,t)=P_{A_\al\ot H/A_\al}(\De_Eb,t)=P_{E\ot H/E}(b,t).
$$
In case $\ga=\la_\al$ we have $(\ga\ot\id)\circ\de=(\io\ot\id)\circ\de_\al$,
whence $(\ga\ot\id)\de a=(\io\ot\id)b$. An application of (P3) for the ring homomorphism
$\io:E\to A_\al$ yields
$$
\la_\al^tP_{A\ot H/A}(\de a,t)=P_{A_\al\ot H/A_\al}\bigl((\io\ot\id)b,t\bigr)
=P_{E\ot H/E}(b,t).
$$

Thus (2.5) is proved. Let $P_{A\ot H/A}(\de a,t)=\sum_{i=0}^nc_it^i$ with
$c_0,\ldots,c_n\in A$. If $\p=\ker\al$ then the sequence
$$
0\mapr{}\p\ot H\mapr{}A\ot H\lmapr4{\al\ot\id}E\ot H
$$
is exact by $K$-projectivity of $H$. Formula (2.5) means that
$\de_\al(c_i)=\la_\al(c_i)$, whence $\de c_i-c_i\ot1\in\ker(\al\ot\id)=\p\ot
H$ for each $i$. This can be rewritten in terms of the $H^*$-module
structure on $A$ as $\xi c_i-\xi(1)c_i\in\p$ for all $\xi\in H^*$, or
$(H^*)^+c_i\sbs\p$. These inclusions hold for every $\p\in\Spec A$ since any prime
ideal $\p$ is the kernel of a homomorphism into a field. Put $J=AV$  where
$V=\sum_{i=0}^n(H^*)^+c_i$. By the above $J$ is contained in the nil radical of $A$.
Since $(H^*)^+$ is an ideal of $H^*$, it is clear that $V$ is an $H^*$-submodule of $A$,
whence so too is $J$. In other words, $J$ is an $H$-costable nil ideal of $A$.
It follows that $J=0$ by the assumption on $A$. Thus $(H^*)^+c_i=0$, and so $c_i\in A^H$.

If $I$ is an $H$-costable ideal of $A$ and $\pi:A\to\A=A/I$ is the canonical
projection then $\de\circ\pi=(\pi\ot\id)\circ\de$, and so
$P_{\A\ot H/\A}(\de\pi a,t)=\pi^tP_{A\ot H/A}(\de a,t)$ for every $a\in A$ by (P3).
This polynomial has all coefficients in $\pi(A^H)\sbs(A/I)^H$. Thus $A/I$ has
invariant characteristic polynomials.

Consider now the general case of Theorem 2.5. Suppose that $\ph$ is given. If
$\ph$ is surjective, then we are done by the previous step. Otherwise consider
the polynomial algebra $A''=A'[X]$ where $X$ is any set of indeterminates.
Extend $\de:A'\to A'\ot H$ to a ring homomorphism $A''\to A''\ot H$ setting
$\de(x)=x\ot1$ for all $x\in X$. This makes $A''$ into an $H$-comodule
commutative algebra. Denote by $R$ the $H$-radical of $A''$. Each $f\in R$ is
nilpotent, whence all coefficients of $f$ are nilpotent \cite{BouA, Ch.~IV,
\S1, Prop.~9}. For every monomial $y$ in $x$'s denote by $R_y\sbs A'$ the
subset consisting of all elements occurring as a coefficient of $y$ in some
$f\in R$. Clearly $R_y$ is a nil ideal of $A'$. Since $R$ is $H$-costable, so
too is $R_y$. It follows that $R_y=0$ since $A'$ is $H$-reduced. As this holds
for every monomial $y$, we conclude that $R=0$. Thus $A''$ is $H$-reduced. We
can extend $\ph$ to a homomorphism of $H$-comodule algebras $\psi:A''\to A$
sending each $x\in X$ to an arbitrary element of $A^H$. Taking the set $X$ big
enough, we can thus obtain a surjective homomorphism $\psi$. As have been
pointed out already, this suffices to complete the proof.
\endproof

\proclaim
Proposition 2.6.
Suppose that $A$ has invariant characteristic polynomials and $m$ is a positive
integer such that $\rk_\q H$ divides $m$ for every $\q\in\Spec K$. If $\a\in(A/I)^H$
where $I$ is an $H$-costable ideal of $A$ then ${m\choose j}\a^j$ belongs to
the image of the canonical map $A^H\to(A/I)^H$ for every $j=0,\ldots,m$.
\endproclaim

\Proof.
Lemma 1.2 gives decompositions $A^H\cong\prod A_i^{H_i}$ and
$(A/I)^H\cong\prod(A_i/I_i)^{H_i}$ where each $I_i$ is an $H_i$-costable ideal
of $A_i$. The assumption on $m$ means that $m$ is divisible by the ranks of the
projective $K_i$-modules $H_i$ for all $i$ occurring here. So it suffices to
prove the assertion of the proposition for the canonical maps
$A_i^{H_i}\to(A_i/I_i)^{H_i}$. We may assume therefore that $H$ has constant rank $n$
and $m=dn$ for some integer $d>0$.

Let $\pi:A\to\A=A/I$ denote the canonical projection, and let $a\in A$ be any
representative of $\a$. By (P3) and (P5)
$$
\pi^tP_{A\ot H/A}(\de a,t)=P_{\A\ot H/\A}(\de\a,t)
=P_{\A\ot H/\A}(\a\ot1,t)=(t-\a)^n.
$$
Thus ${n\choose l}(-\a)^l=\pi(c_{n-l})$ for $l=0,\ldots,n$ where
$\sum_{l=0}^nc_lt^l=P_{A\ot H/A}(\de a,t)$. By the hypothesis
$c_0,\ldots,c_n\in A^H$. Now
$$
\textsum\limits_{j=0}^m{m\choose j}\a^jt^{m-j}=\textstyle(t+\a)^m=(t+\a)^{dn}
=\bigl(\,\sum\limits_{l=0}^n{n\choose l}\a^lt^{n-l}\,\bigr)^d,
$$
and therefore the elements ${m\choose j}\a^j$ with $j=0,\ldots,m$ belong to
the subalgebra of $A/I$ generated by the images of $c_0,\ldots,c_n$.
\endproof

\Remark.
In particular, $\a^m$ is in the image of $A^H\to(A/I)^H$ for every
$\a\in(A/I)^H$. This can be regarded as a {\it geometric reductivity} of finite
Hopf algebras. Classically, the geometric reductivity is a property of a
reductive algebraic group $G$ over a field of characteristic $p>0$ which serves
as a substitute for linear reductivity \cite{Mum65}, \cite{Hab75}. One of its
formulations is as follows: if $G$ operates rationally on a commutative algebra
$A$ as a group of automorphisms and $I$ is a $G$-stable ideal of $A$, then for
every $G$-invariant $\a\in(A/I)^G$ there exists an integer $m>0$ such that
$\a^m$ has a representative in $A^G$. The geometric reductivity of finite group
schemes is a byproduct of the norm \cite{Wat}.
\endremark

\proclaim
Proposition 2.7.
$A$ is integral over $A^H$ in any of the following two cases{\rm:}

\item(a)
$A$ has invariant characteristic polynomials,

\item(b)
the $H$-radical $R$ of $A$ is $\Z$-torsion.

\endproclaim

\Proof.
(a) It suffices to show that each $A_i$ in Lemma 1.2 is integral over $A_i^{H_i}$
by \cite{Bou, Ch.~V, \S1, Prop.~3}. So we may assume without loss of generality that
$H$ is $K$-projective of constant rank $n$. Given $a\in A$, let
$P_{A\ot H/A}(\de a,t)=\sum_{i=0}^nc_it^i$. By the hypothesis
$c_0,\ldots,c_n\in A^H$. Property (P1) ensures that $\sum_{i=0}^n(c_i\ot1)(\de a)^i=0$
in $A\ot H$. Applying the ring homomorphism $\id\ot\ep:A\ot H\to A$ to both sides of
this equality, we get $\sum_{i=0}^nc_ia^i=0$ since $(\id\ot\ep)\circ\de=\id$. Here
$c_n=1$, and so $a$ is integral over $A^H$.

(b) The $H$-comodule algebra $A/R$ is $H$-reduced, and so $A/R$ has invariant
characteristic polynomials by Theorem 2.5. Thus $A/R$ is integral over $(A/R)^H$ by
case (a). Let $B\subset A$ denote the preimage of $(A/R)^H$ with respect to the
canonical map $A\to A/R$. Clearly $B$ is a subalgebra of $A$ containing $A^H$ and $A$
is integral over $B$. It remains to prove that $B$ is integral over $A^H$
\cite{Bou, Ch.~V, \S1, Prop.~6}. Let $b\in B$. In view of the exact sequence
$$
R\ot H\to A\ot H\to(A/R)\ot H\to0
$$
the element $\de b-b\ot1\in A\ot H$ which goes to 0 in $(A/R)\ot H$ belongs
to the image of $R\ot H$. Thus $\de b=b\ot1+x$ where
$x=x_1\ot h_1+\ldots+x_l\ot h_l$ with $x_1,\ldots,x_l\in R$ and
$h_1,\ldots,h_l\in H$. Since $R$ is a nil ideal, the ideal of $A$ generated
by finitely many $x_1,\ldots,x_l$ is nilpotent. There exists therefore an
integer $n>0$ such that $x^n=0$. By the hypothesis there exists an
integer $e>0$ such that $ex_j=0$ for all $j$, hence also $ex=0$. Take an integer
$m>0$ such that $e$ divides all binomial coefficients $m\choose i$ for
$i=1,\ldots,n$ (for instance, $m=e\cdot n!$ will do). Noting that $b\ot1$ is
central in the algebra $A\ot H$, we obtain
$\de(b^m)=(\de b)^m=\sum_{i=0}^m{m\choose i}(b^{m-i}\ot1)x^i=b^m\ot1$. This
shows that $b^m\in A^H$. Hence $b$ is integral over the subring $A^H$.
\endproof

\Remark.
If $A$ is an algebra over a field of characteristic $p>0$ then condition (b)
is fulfilled since $pA=0$.
\endremark

\section
3. The orbital subalgebras and the quotient map

For each $\p\in\Spec A$ denote by $k(\p)$ the residue field of the local ring
$A_\p$. The same notation will be used for commutative rings other than $A$.
Let $\al_\p:A\to k(\p)$ be the canonical ring homomorphism. Denote by $O(\p)$ the commutative
right coideal subalgebra $A_{\al_\p}$ of the Hopf algebra $k(\p)\ot H$ over $k(\p)$ defined
in (2.1). We call $O(\p)$ the {\it orbital subalgebra} associated with $\p$. In geometric terms
$\Spec O(\p)$ is isomorphic with a closed subscheme of $\Spec k(\p)\times_{\Spec K}\Spec A$; when
$H$ is commutative and $G=\Spec H$, this subscheme is called the $G$-orbit of $\p$. Let
$\de_\p=\de_{\al_\p}$ as in (2.2). Thus $\de_\p$ is a ring homomorphism $A\to
k(\p)\ot H$ and $O(\p)=(k(\p)\ot1)\cdot\de_\p(A)$.

\proclaim
Lemma 3.1.
The kernel of $\de_\p$ is the largest $H$-costable ideal of $A$ contained in $\p$.
\endproclaim

\Proof.
In view of Lemma 2.1 $\de_\p$ is a homomorphism of $H$-comodule algebras. Hence $\ker\de_\p$
is an $H$-costable ideal of $A$. The composite of two arrows at the top of commutative diagram
$$
\vcenter{\diagram{
A&\lmapr1\de&A\ot H&\lmapr5{\al_\p\ot\id}&k(\p)\ot H\cr
\diagramskip
&\tsear{-2pt}\id\hidewidth&\mapd{}{\id\ot\ep}&&\mapd{}{\id\ot\ep}\cr
\diagramskip
&&A&\hidewidth\lmapr8{\al_\p}\hidewidth&k(\p)\cr
}}\eqno(3.1)
$$
equals $\de_\p$, whence $(\id\ot\ep)\circ\de_\p=\al_\p$. This yields
$\ker\de_\p\sbs\ker\al_\p=\p$. If $I$ is any $H$-costable ideal of
$A$ and $I\sbs\p$, then $\de_\p(I)\sbs(\al_\p\ot\id)(I\ot H)=0$.
\endproof

\proclaim
Lemma 3.2.
Let $\ga:A\to B$ be a homomorphism of commutative $H$-comodule algebras, and
let $\p=\ga^{-1}(\q)$ where $\q\in\Spec B$. Then:

\item(i)
$(\g\ot\id)\circ\de_\p=\de_\q\circ\ga$ where $\g:k(\p)\to k(\q)$ is induced by
$\ga$.

\item(ii)
$k(\q)\ot_{k(\p)}O(\p)$ is canonically embedded into $O(\q)$ as a
$k(\q)$-subalgebra.

\item(iii)
If $B=\ga(A)B^H$ then $O(\q)\cong k(\q)\ot_{k(\p)}O(\p)$.

\item(iv)
If $\ga a_1,\ldots,\ga a_d$ generate $B$ as a $B^H$-module for some
$a_1,\ldots,a_d\in A$, then $\de_\p a_1,\ldots,\de_\p a_d$ span $O(\p)$ over
$k(\p)$.

\endproclaim

\Proof.
Part (i) is seen from the commutative diagram below in which the composite of
maps at the top equals $\de_\p$ and that at the bottom equals $\de_\q$:
$$
\diagram{
A&\lmapr3\de&A\ot H&\lmapr6{\al_\p\ot\id}&k(\p)\ot H\cr
\diagramskip
\mapd{}\ga&&\mapd{}{\ga\ot\id}&&\mapd{}{\g\ot\id}\cr
\diagramskip
B&\lmapr3\de&B\ot H&\lmapr6{\al_\q\ot\id}&k(\q)\ot H.\cr
}
$$

(ii) The inclusion of $O(\p)$ in $k(\p)\ot H$ induces an embedding of
$k(\q)$-algebras
$$
\ph:k(\q)\ot_{k(\p)}O(\p)\to k(\q)\ot_{k(\p)}\bigl(k(\p)\ot H\bigr)\cong
k(\q)\ot H.
$$
If $a\in A$, then $\ph(1\ot\de_\p a)=(\g\ot\id)(\de_\p a)=\de_\q(\ga a)$. This
shows that $\im\ph\sbs O(\q)$.

(iii) By the computation in (ii) $\de_\q(\ga a)\in\im\ph$ for all $a\in A$. If
$b\in B^H$ then $ \de_\q b=\al_\q(b)\ot1\in\im\ph$ as well. Hence
$\de_\q(B)\sbs\im\ph$ provided that $B=\ga(A)B^H$. By the definition
$O(\q)=(k(\q)\ot1)\cdot\de_\q(B)\sbs\im\ph$.

(iv) The hypothesis of (iii) is fulfilled, so that $\ph$ is an isomorphism.
Repeating the arguments in (iii), we obtain, more precisely, that $O(\q)$ is
spanned over $k(\q)$ by $\de_\q(\ga a_1),\ldots,\de_\q(\ga a_d)$ which are the
images of $1\ot\de_\p a_1,\ldots,1\ot\de_\p a_d$ under $\ph$. Thus
$O(\q)=\ph\bigl(k(\q)\ot_{k(\p)}\!V\bigr)$ where $V\sbs O(\p)$ is the
$k(\p)$-linear span of $\de_\p a_1,\ldots,\de_\p a_d$, and so $V=O(\p)$.
\endproof

\Remark.
If $B=\ga(A)B^H$ and $\g$ is an isomorphism $k(\p)\cong k(\q)$, then
$O(\p)\cong O(\q)$. This is the case when $B$ is either $A/I$ with $I$ an
$H$-costable ideal of $A$ or $B=S^{-1}A$ with $S$ a multiplicatively closed
subset of $A^H$ and $\ga$ is the canonical homomorphism in both cases. The
same is true when $\ga$ is a composite of homomorphisms of these two types.
\endremark

The embedding $A^H\to A$ induces a map $\pi:\Spec A\to\Spec A^H$ defined by
the rule $\p\mapsto\p\cap A^H$ for each prime ideal $\p$ of $A$. In the case
when $H$ is commutative, $\Spec A^H$ coincides with the quotient
$\Spec A/\Spec H$ considered in \cite{De}, \cite{Mum70}. The known results
describing the properties of $\pi$ can be extended to noncommutative Hopf
algebras $H$.

In the rest of this section we assume that {\it$A$ has invariant
characteristic polynomials}. Note that the integrality of $A$ over $A^H$
(Proposition 2.7) implies that $\pi$ is surjective and {\it closed}, that is,
$\pi$ maps closed subsets of $\Spec A$ onto closed subsets of $\Spec A^H$
\cite{Bou, Ch.~V, \S2, Th.~1, Remark 2}. A less obvious property is the {\it
openness} of $\pi$ (the open subsets of $\Spec A$ are mapped onto open subsets
of $\Spec A^H$).

\proclaim
Theorem 3.3.
Let $\p\in\Spec A$. Then:

\item(i)
There are only finitely many $\p'\in\Spec A$ such that
$\p'\cap A^H=\p\cap A^H$. These are precisely the ideals
$\p'=\de_\p^{-1}(\m)$ for some maximal ideal $\m$ of $O(\p)$.

\item(ii)
The map $\pi:\Spec A\to\Spec A^H$ is open.

\item(iii)
If $\q'\in\Spec A^H$ satisfies $\q'\subset\p\cap A^H$
then $\q'=\p'\cap A^H$ for some $\p'\in\Spec A$ such that $\p'\subset\p$
{\rm(}Going-down property{\rm)}.

\endproclaim

\Proof.
In view of Lemma 1.2 the proof is reduced to the case where $H$ is
$K$-projective of constant rank $n$.

(i) Put $\q=\p\cap A^H$. Suppose that $\p'=\de_\p^{-1}(\m)$ where $\m$ is a maximal
ideal of $O(\p)$. Then $\p'$ is a prime ideal of $A$. Given $a\in A^H$, we have
$\de_\p(a)=\al_\p(a)\ot1$. Since $k(\p)\ot1$ is a subfield in $O(\p)$, it has
zero intersection with $\m$. Hence $\de_\p(a)\in\m$ if and only if $\al_\p(a)=0$,
that is, $a\in\p'$ if and only if $a\in\p$. This shows that $\p'\cap A^H=\q$.

Suppose now that $\p'$ is any prime ideal of $A$ such that $\p'\cap A^H=\q$.
Let $\m_1,\ldots,\m_s$ be all maximal ideals of the finite dimensional algebra $O(\p)$,
and put $\p_j=\de_\p^{-1}(\m_j)$. Suppose that $\p'\not\subset\p_j$ for every $j$.
Then $\p'\not\subset\p_1\cup\ldots\cup\p_s$ by \cite{Bou, Ch.~II, \S1, Prop.~1}. Take
any $a\in\p'$ outside of every $\p_j$ and consider
$P_{A\ot H/A}(\de a,t)=\sum_{i=0}^nc_it^{n-i}$. By the assumption $c_i\in A^H$
for all $i$. Note that $c_n=(-1)^nN_{A\ot H/A}(\de a)$. As
$\sum_{i=0}^nc_ia^{n-i}=0$ (case (a) of Proposition 2.7), we have $c_n\in
aA\subset\p'$. Hence $c_n\in\q$, and $\al_\p(c_n)=0$. On the other hand
$\al_\p(c_n)=(-1)^nN_{k(\p)\ot H/k(\p)}(\de_\p a)$ by (P3). By the choice of $a$ the
element $\de_\p a$ belongs to none of $\m_1,\ldots,\m_s$ and therefore is
invertible in $O(\p)$. Then $\de_\p a$ is invertible also in $k(\p)\ot H$,
and (P7) shows that $\al_\p(c_n)$ is invertible in $k(\p)$, yielding a
contradiction. Thus $\p'\subset\p_j$ for at least one $j$. Now $\p_j\cap
A^H=\q$ by the previous step in the proof, and the integrality of $A$ over
$A^H$ implies that $\p'=\p_j$ \cite{Bou, Ch.~V, \S2, Cor.~1 to Lemma 2}.

(ii) For $a\in A$ put $D_a=\{\p\in\Spec A\mid a\notin\p\}$. The subsets $D_a$
form a base of topology on $\Spec A$. So it suffices to prove that $\pi(D_a)$
is an open subset of $\Spec A^H$ for each $a$. Let
$P_{A\ot H/A}(\de a,t)=\sum_{i=0}^nc_it^{n-i}$, so that $c_i\in A^H$ for
all $i$. Denote by $I_a$ the ideal of $A^H$ generated by $c_1,\ldots,c_n$. We will
show that
$$
\pi(D_a)=\{\q\in\Spec A^H\mid I_a\not\subset\q\}.\eqno(3.2)
$$

Suppose that $\p\in\Spec A$. If $c_1,\ldots,c_n\in\p$ then
$a^n=-\sum_{i=1}^nc_ia^{n-i}\in\p$, whence $a\in\p$. Therefore $a\notin\p$
implies $I_a\not\subset\p\cap A^H$.

Suppose that $\q\in\Spec A^H$. Let $\p\in\Spec A$ be any prime ideal lying
above $\q$. Denote by $\m_1,\ldots,\m_s$ all maximal ideals of $O(\p)$, and put
$\p_j=\de_\p^{-1}(\m_j)$. By (i) $\p_j\cap A^H=\q$ for each $j$. Suppose that
$a\in\p_j$ for all $j$. Then $\de_\p a\in\m_j$ for all $j$, whence $\de_\p a$
is a nilpotent element of $O(\p)$. Property (P6) yields
$P_{k(\p)\ot H/k(\p)}(\de_\p a)=t^n$. On the other hand
$P_{k(\p)\ot H/k(\p)}(\de_\p a,t)=\al_\p^tP_{A\ot H/A}(\de a,t)$ by (P3), and
it follows that $\al_\p(c_i)=0$, that is, $c_i\in\p$ for each $i=1,\ldots,n$.
Thus $I_a\not\subset\q$ implies that $a\notin\p_j$ for at least one $j$.
This proves (3.2), and the openness of $\pi(D_a)$ is clear.

Part (iii) is an easy consequence of (i) and (ii). Denote by
$\p_1,\ldots,\p_r\in\Spec A$ all prime ideals lying above $\q'$. For each
$j=1,\ldots,r$ consider the open subset $U_j\subset\Spec A$ consisting of
those prime ideals of $A$ which do not contain $\p_j$. The intersection
$U=U_1\cap\cdots\cap U_r$ is also open in $\Spec A$. Now $\pi(U)$ is open in
$\Spec A^H$ by (ii), and $\q'\notin\pi(U)$ by the choice of $\p_1,\ldots,\p_r$.
Put $\q=\p\cap A^H$. Since $\q'\subset\q$, we get $\q\notin\pi(U)$ as well.
Hence $\p\notin U$. It follows that $\p_j\subset\p$ for at least one $j$, and
we may take $\p'=\p_j$.
\endproof

For $\p\in\Spec A$ denote by $\dim O(\p)$ the dimension of $O(\p)$ as a vector
space over $k(\p)$. We say that $\p$ is {\it$H$-regular} if
$\dim O(\p')=\dim O(\p)$ for all $\p'$ in a suitable neighbourhood of $\p$ in
$\Spec A$.

\proclaim
Lemma 3.4.
Suppose that $\de_\p a_1,\ldots,\de_\p a_d$ are a basis for $O(\p)$ over
$k(\p)$ where $a_1,\ldots,a_d\in A$. Then{\rm:}

\item(i)
If $\p'\cap A^H=\p\cap A^H$ for $\p'\in\Spec A$ then
$\de_{\p'}a_1,\ldots,\de_{\p'}a_d$ are a basis for $O(\p')$ over $k(\p')$.

\item(ii)
There exists a neighbourhood $U$ of $\p$ in $\Spec A$ such that
$\de_{\p'}a_1,\ldots,\de_{\p'}a_d$ are linearly independent over $k(\p')$
for each $\p'\in U$.

\item(iii)
If $\p$ is $H$-regular, then there exists $s\in A^H\setm\p$ such that
$\de_{\p'}a_1,\ldots,\de_{\p'}a_d$ are a basis for $O(\p')$ over $k(\p')$
whenever $\p'\in\Spec A$ and $s\notin\p'$.

\endproclaim

\Proof.
(i) By Theorem 3.3 $\p'=\de_\p^{-1}(\m)$ for some maximal ideal $\m$ of
$O(\p)$. We apply Lemma 3.2 taking $B=O(\p)$, $\ga=\de_\p$, $\q=\m$, and $\p'$
in place of $\p$. Note that $k(\p)\ot1\sbs B^H$. Hence the hypothesis in (iv)
of Lemma 3.2 is fulfilled, and we deduce that
$\de_{\p'}a_1,\ldots,\de_{\p'}a_d$ span $O(\p')$ over $k(\p')$. It follows
that $\dim O(\p')\le\dim O(\p)$. By symmetry between $\p$ and $\p'$ there
is also the opposite inequality, which gives $\dim O(\p')=d$.

(ii) There exist $k(\p)$-linear functions $\xi_1,\ldots,\xi_d:k(\p)\ot H\to
k(\p)$ such that $\xi_i(\de_\p a_j)$ is 1 when $i=j$ and 0 otherwise. Since
$H$ is a finitely projective $K$-module, we have
$$
\Hom_{k(\p)}\bigl(k(\p)\ot H,k(\p)\bigr)\cong\Hom_K\bigl(H,k(\p)\bigr)\cong k(\p)\ot H^*.
$$
We can find therefore $\eta_1,\ldots,\eta_d\in H^*$ such that the $d\times d$
matrix with entries $(\id\ot\eta_i)(\de_\p a_j)\in k(\p)$ is invertible.
Denote by $M$ and $M(\p')$, where $\p'\in\Spec A$, the $d\times d$ matrices
with entries $(\id\ot\eta_i)(\de a_j)\in A$ and
$(\id\ot\eta_i)(\de_{\p'}a_j)\in k(\p')$, respectively. Thus $M(\p')$ is
obtained by applying $\al_{\p'}$ to all entries of $M$. Put $u=\det M\in A$.
Then $\al_{\p'}(u)=\det M(\p')$, and so $M(\p')$ is invertible if and only if
$u\notin\p'$. In particular, $u\notin\p$. If $u\notin\p'$ then
$\de_{\p'}a_1,\ldots,\de_{\p'}a_d$ are linearly independent over $k(\p')$. The
set $D_u=\{\p'\in\Spec A\mid u\notin\p'\}$ is the required open neighbourhood
of $\p$ in $\Spec A$.

(iii) Let $U$ be as in (ii). Taking a smaller neighbourhood, if necessary, we
may assume that $\dim O(\p')=\dim O(\p)$ for all $\p'\in U$. Then
$\de_{\p'}a_1,\ldots,\de_{\p'}a_d$ are a basis for $O(\p')$ over $k(\p')$
whenever $\p'\in U$. By Theorem 3.3 $\pi(U)$ is a neighbourhood of
$\q=\p\cap A^H$ in $\Spec A^H$. There exists $s\in A^H$ such that $\q\in
D_s\sbs\pi(U)$ where $D_s$ is the subset of prime ideals of $A^H$ not
containing $s$. If $\p'\in\Spec A$ and $s\notin\p'$ then
$\p'\cap A^H=\p''\cap A^H$ for some $\p''\in U$. Applying (i) with $\p''$ in
place of $\p$, we get the desired conclusion.
\endproof

\proclaim
Proposition 3.5.
{\rm(i)}
The function $\p\mapsto\dim O(\p)$ is lower semicontinuous, that is, for each
integer $m$ the subset $\{\p\in\Spec A\mid\dim O(\p)\ge m\}$ is open in $\Spec
A$.

\item(ii)
The subset of all $H$-regular prime ideals is dense and open in $\Spec A$.

\item(iii)
If $\p'\cap A^H=\p\cap A^H$ for $\p,\p'\in\Spec A$ then
$\dim O(\p')=\dim O(\p)$. For $\p'$ to be $H$-regular, it is necessary and
sufficient that so be $\p$.

\endproclaim

\Proof.
(i) Clear from Lemma 3.4(ii).

(ii) The openness property is clear from the definition. Suppose that
$U\subset\Spec A$ is any nonempty open subset. Put $m=\max_{\p\in U}\dim O(\p)$.
This maximum makes sense since $\dim O(\p)\le n$ for every $\p$ where $n$ is
the number of generators for $H$ as a $K$-module. The subset
$V=\{\p\in\Spec A\mid\dim O(\p)\ge m\}$ is open by (i). Now $U\cap V$ is
nonempty and consists of $H$-regular primes. This verifies the density
property.

(iii) The equality of dimensions is proved in Lemma 3.4(i). If $\p$ is
$H$-regular and $s$ is as in Lemma 3.4(iii), then $U=\{\p''\in\Spec A\mid
s\notin\p''\}$ is a neighbourhood of $\p'$ in $\Spec A$ such that the algebras
$O(\p'')$ have the same dimension for all $\p''\in U$. In this case
$\p'$ is $H$-regular.
\endproof

\proclaim
Lemma 3.6.
Suppose that $A^H$ is local with a maximal ideal $\q$. If there exists an
$H$-regular $\p\in\Spec A$ lying above $\q$ then all prime ideals of $A$ are
$H$-regular.
\endproclaim

\Proof.
If $\p'$ is any maximal ideal of $A$ then $\p'\cap A^H=\q$ since $A$ is integral over
$A^H$. Then all maximal ideals of $A$ are $H$-regular by Proposition 3.5(iii),
and so are all prime ideals by (ii).
\endproof

\section
4. Projectivity over the subring of invariants

We say that $A$ is {\it$H$-simple} if $A$ has no nonzero proper $H$-costable ideals.

\proclaim
Lemma 4.1.
Suppose that $A$ is $H$-simple. Then $A^H$ is a field and $A$ is finite
dimensional over $A^H$. If $\p\in\Spec A$ and $a_1,\ldots,a_d\in A$ are such
that $\de_\p a_1,\ldots,\de_\p a_d$ form a basis for $O(\p)$ over $k(\p)$
then $a_1,\ldots,a_d$ are a basis for $A$ over $A^H$.
\endproclaim

\Proof.
The hypothesis means that $A$ is a simple object of $\M_A^H$. Every $A$-module
endomorphism of $A$ is obtained by the rule $a\mapsto ca$ for some $c\in A$.
In order that this map $A\to A$ commute with the coaction of $H$, it is
necessary and sufficient that $c\in A^H$. In other words $A^H$ is identified
with the endomorphism ring of $A$ as an object of $\M_A^H$. By Schur's lemma
$A^H$ is a field. Given $\p$, the field $A^H$ is embedded in $k(\p)$. We may
regard $A'=k(\p)\ot_{A^H}A$ as an object of $\M_A^H$ using the operations on
the second tensorand. Thus $k(\p)\ot1\sbs A'^H$. Taking any basis of $k(\p)$
over $A^H$, we present $A'$ as a direct sum of copies of $A$. Hence $A'$ is a
semisimple object of $\M_A^H$, and every its subobject is a sum of simple
subobjects isomorphic to $A$. Now any morphism $A\to A'$ in $\M_A^H$ is of the
form $a\mapsto c\ot a$ for some $c\in k(\p)$. So it follows that every
subobject of $A'$ is of the form $V\ot_{A^H}A$ where $V\sbs k(\p)$ is an
$A^H$-subspace. Denote by $\ph:A'\to k(\p)\ot H$ the $k(\p)$-linear map
extending $\de_\p:A\to k(\p)\ot H$. Then $\ph$ is a homomorphism of
$H$-comodule algebras, and so $I=\ker\ph$ is an $H$-costable ideal of $A'$.
Hence $I$ is a subobject of $A'$ in $\M_A^H$. We deduce that $I=V\ot_{A^H}A$
for some $V$ as above. On the other hand, $\ph$ is injective on $k(\p)\ot1$.
Since $V\ot1\sbs I$, we get $V=0$. Thus $\ph$ is injective. The image of $\ph$
coincides with $O(\p)$ by the definition.  If $\de_\p a_1,\ldots,\de_\p a_d$
form a basis for $O(\p)$ over $k(\p)$ then $1\ot a_1,\ldots,1\ot a_d$ form a
basis for $A'$ over $k(\p)$, and the final assertion is clear.
\endproof

\proclaim
Lemma 4.2.
Suppose that $A$ has a maximal ideal $\p$ which contains no nonzero
$H$-costable ideals of $A$. Then $A$ is $H$-simple.
\endproclaim

\Proof.
The $H$-radical of $A$ is contained in every prime ideal of $A$. The hypothesis
of the lemma implies therefore that $A$ is $H$-reduced. By Theorem 2.5 $A$ is
integral over $A^H$. Next, the ideal of $A$ generated by $\p\cap A^H$ is
$H$-costable and is contained in $\p$. Hence $\p\cap A^H=0$. The maximality of
$\p$ ensures that $(0)$ is a maximal ideal of $A^H$ \cite{Bou, Ch.~V, \S2,
Prop.~1}. In other words $A^H$ is a field. Suppose that $I$ is any proper
$H$-costable ideals of $A$. Then $I\sbs\p'$ for some maximal ideal $\p'$ of $A$.
By Lemma 3.1 $\de_{\p'}(I)=0$. As $\p'\cap A^H=0$, Theorem 3.3 shows that
$\p=\de_{\p'}^{-1}(\m)$ for some maximal ideal $\m$ of $O(\p')$. But then
$I\sbs\p$, and so $I=0$ by the hypothesis.
\endproof

Denote by $\M'$ the full subcategory of $\M_A^H$ consisting of right
$(H,A)$-Hopf modules $M$ such that $M=M^HA$. The next result generalizes
\cite{Sk02, Th.~2.1, Prop.~3.2} where $K$ was supposed to be an algebraically
closed field, $A$ a finitely generated integral domain and $H$ a commutative
Hopf algebra.

\proclaim
Theorem 4.3.
Suppose that $A$ is $H$-reduced and the function $\p\mapsto\dim O(\p)$ is locally constant on the
whole $\Spec A$. Then{\rm:}

\item(i)
$A$ is a finitely generated projective $A^H$-module whose rank at $\q\in\Spec A^H$ is
equal to $\dim O(\p)$ where $\p$ is any prime ideal of $A$ lying above $\q$.

\item(ii)
The functor $M\mapsto M^H$ is an equivalence between $\M'$ and the category of
$A^H$-modules. The inverse functor is $N\mapsto N\otimes_{A^H}A$.

\item(iii)
The assignment $I\mapsto I\cap A^H$ establishes a bijection between the
$H$-costable ideals of $A$ and the ideals of $A^H$. The inverse
correspondence is $J\mapsto JA$.

\endproclaim

\Proof.
For every $s\in A^H$ the localization $A_s$ is an $H$-reduced $H$-comodule
algebra. Given $M\in\M_A^H$, one has $M_s\in\M_{A_s}^H$ by Lemma 1.1. If,
moreover, $M=M^HA$ then $M_s=M_s^HA$. To prove (i) it suffices, by \cite{Bou,
Ch.~II, \S5, Th.~1}, to show that for every $\q\in\Spec A^H$ there exists
$s\in A^H$ such that $s\notin\q$ and $A_s$ is a free $A_s^H$-module of rank
$d=\dim O(\p)$ (by Proposition 3.5 $d$ does not depend on a choice of $\p$
above $\q$). In (ii) one first notes that the two functors are well defined;
the $H$-comodule structure on $N\otimes_{A^H}A$ is given by means of the map
$$
\id\otimes\,\de:N\otimes_{A^H}A\to N\otimes_{A^H}(A\otimes H)\cong
(N\otimes_{A^H}A)\otimes H.
$$
For every $M\in\M'$ and every $A^H$-module $N$ define
$$
\Psi_M:M^H\otimes_{A^H}A\to M,\qquad\qquad\Phi_N:N\to(N\otimes_{A^H}A)^H.
$$
by $v\ot a\mapsto va$ and $v\mapsto v\ot1$, respectively. Then $\Psi_M$ is a
morphism in $\M'$ and $\Phi_N$ is a homomorphism of $A^H$-modules. The map
$\Psi_M$ is bijective if and only if for every $\q\in\Spec A^H$ there exists
$s\in A^H$ such that $s\notin\q$ and $\Psi_M\ot_AA_s$ is bijective (cf.
\cite{Bou, Ch.~II, \S3, Th.~1}). Note that $\Psi_M\ot_AA_s$ can be identified
with the map $M_s^H\otimes_{A_s^H}A_s\to M_s$ such that $v\ot a\mapsto va$.
The bijectivity of $\Phi_N$ can be verified similarly. Thus, in proving (i)
and (ii), we may pass to suitable localizations $A_s$.

By the hypothesis all prime ideals of $A$ are $H$-regular. Taking $s$ as in
Lemma 3.4(iii) with respect to any chosen prime ideal of $A$ and replacing
$A$ with $A_s$, we reduce the proof to the case in which there exist
$a_1,\ldots,a_d\in A$ such that for every $\p\in\Spec A$ the elements
$\de_\p a_1,\ldots,\de_\p a_d$ form a basis for $O(\p)$ over $k(\p)$. We fix
$a_1,\ldots,a_d$.

Let us regard $A\ot H$ as an $A$-module via left multiplications on the
first tensorand. This module is finitely projective since so is the $K$-module $H$.
Consider its submodule $E\sbs A\ot H$ generated by $\de a_1,\ldots,\de a_d$. If
$\p\in\Spec A$ then the $A_\p$-module
$A_\p\ot H$ obtained by localizing $A\ot H$ at $\p$ is free of finite rank.
Furthermore, $E_\p$ may be identified with an $A_\p$-submodule of $A_\p\ot
H$ generated by $d$ elements whose images in $k(\p)\ot H$ are linearly
independent over $k(\p)$. Then $E_\p$ is an $A_\p$-module direct summand of
$A_\p\ot H$ \cite{Bou, Ch.~II, \S3, Cor.~1 to Prop.~5}. In particular, $E_\p$
is a free $A_\p$-module of finite rank. If $\ph:A^d\to E$ is the $A$-module
homomorphism sending the standard generators of the rank $d$ free $A$-module
$A^d$ to $\de a_1,\ldots,\de a_d$ then the localizations of $\ph$ at prime
ideals of $A$ are all isomorphisms. Hence $\ph$ is itself an isomorphism,
i.e., $E$ is a free $A$-module with $\de a_1,\ldots,\de a_d$ as its basis. By
\cite{Bou, Ch. II, \S3, Cor.~1 to Prop. 12} $E$ is an $A$-module direct
summand of $A\ot H$.

\proclaim
Claim 1.
Let $I(\p)$ be the largest $H$-costable ideal of $A$ contained in a prime ideal
$\p$. Then $\de(A)\sbs E+I(\p)\ot H$.
\endproclaim

Put $I=I(\p)$ and $\q=\p\cap A^H$. The prime ideal $\p A_\q$ of $A_\q$ lies
above the maximal ideal $\q A^H_\q$ of $A^H_\q$. By Theorem 2.5 $A_\q$ is
integral over $A^H_\q$, whence $\p A_\q$ is a maximal ideal of $A_\q$. Next,
$IA_\q$ is the largest $H$-costable ideal of $A_\q$ contained in $\p A_\q$. Lemma
4.2 shows that the $H$-comodule algebra $A'=A_\q/IA_\q$ is $H$-simple. If
$\p'=\p A_\q/IA_\q$, then there is an isomorphism $O(\p)\to O(\p')$ (see Lemma
3.2) such that $\de_\p a_i\mapsto\de_{\p'}a'_i$ where $a'_i$ denotes the
image of $a_i$ in $A'$ for each $i=1,\ldots,d$. Hence
$\de_{\p'}a'_1,\ldots,\de_{\p'}a'_d$ form a basis for $O(\p')$ over $k(\p')$.
Applying Lemma 4.1, we deduce that $A'$ is spanned over $A'^H$ by
$a'_1,\ldots,a'_d$. As $\de(A'^H)\sbs A'\ot1$, it follows that $\de(A')$ is
contained in the $A'$-submodule $E'\sbs A'\ot H$ generated by $\de
a'_1,\ldots,\de a'_d$. Clearly $E'$ coincides with the image of the
localization $E_\q$ under canonical map $A_\q\ot H\to A'\ot H$. The kernel of
the latter map is $IA_\q\ot H$. Hence
$$
\de(A_\q)\sbs E_\q+IA_\q\ot H.\eqno(4.1)
$$
Put $J=\{b\in A\mid ub\in I$ for some $u\in A^H\setm\q\}$. For every
projective $A$-module $P$ we have
$$
JP=\{x\in P\mid ux\in IP{\fam0\ for\ some\ }u\in A^H\setm\q\}\eqno(4.2)
$$
(this is clear for free $A$-modules, hence also for their direct summands). If
$b\in J$ and $ub\in I$ for $u\in A^H\setm\q$ then $(u\ot1)\de(b)=\de(ub)\in I\ot H$.
Equality (4.2) applied to $P=A\ot H$ yields $\de(b)\in J\ot H$. Thus $J$ is
an $H$-costable ideal of $A$. On the other hand $J\sbs\p$ since $I\sbs\p$ and
$\p$ is prime. It follows that $J\sbs I$ by the maximality of $I$. In fact
$J=I$ since the opposite inclusion is obvious. Let $a\in A$. By (4.1) there
exists $u\in A^H\setm\q$ such that $(u\ot1)\de(a)\in E+I\ot H$. Taking
$P=(A\ot H)/E$, we deduce from (4.2) that $\de(a)\in E+J\ot H=E+I\ot H$. Claim
1 is proved.

\proclaim
Claim 2.
Put $B=(A\ot1)\cdot\de(A)\sbs A\ot H$. Then $B$ is an $A$-module direct
summand of $A\ot H$ and is a free $A$-module with $\de a_1,\ldots,\de a_d$ as its basis.
\endproclaim

Put $J=\bigcap_{\p\in\Spec A}I(\p)$. Then
$$
JP=\textstyle\bigcap\limits_{\p\in\Spec A}I(\p)P\eqno(4.3)
$$
for every projective $A$-module $P$ (the verification reduces again to free
$A$-modules). Taking $P=A\ot H$, we see that $\de(J)\sbs\bigcap_{\p\in\Spec
A}I(\p)\ot H=J\ot H$. Thus $J$ is an $H$-costable ideal of $A$. Since
$J\sbs\p$ for every $\p\in\Spec A$, all elements of $J$ are nilpotent. Hence
$J=0$ by the assumption on $A$. Using Claim 1 and (4.3) with $P=(A\ot H)/E$,
we get
$$
\de(A)\sbs\textstyle\bigcap\limits_{\p\in\Spec A}(E+I(\p)\ot H)=E+J\ot H=E.
$$
Hence $B\sbs E$ as well. Then $B=E$ since the opposite inclusion is obvious.
This completes the proof of Claim 2.

For each $A$-module $M$ consider $M\ot H$ as an $A\ot H$-module by means of
the operation $(v\ot g)(a\ot h)=va\ot gh$ where $v\in M$, $a\in A$ and $g,h\in
H$. Put $E_M=(M\ot1)\cdot\de(A)\sbs M\ot H$, which is a $B$-submodule. Clearly
$E_M$ coincides with the image of the map
$$
M\ot_AB\to M\ot_A(A\ot H)\cong M\ot H
$$
induced by the inclusion of $B$ into $A\ot H$. Since $B$ is an $A$-module direct
summand of $A\ot H$, the map above is injective, and so $E_M\cong M\ot_AB$.

\proclaim
Claim 3.
Each element of $E_M$ can be written as $\sum_{i=1}^d(v_i\ot1)\cdot\de(a_i)$
with uniquely determined $v_1,\ldots,v_d\in M$.
\endproclaim

This is immediate from the freeness of $B$ over $A$.

\proclaim
Claim 4.
If $M\in\M'$ then each element of $M$ can be written as $\sum_{i=1}^dv_ia_i$ with uniquely
determined $v_1,\ldots,v_d\in M^H$.
\endproclaim

If $v_1,\ldots,v_d\in M^H$ and $\sum v_ia_i=0$ then
$\sum(v_i\ot1)\de(a_i)=\de(\sum v_ia_i)=0$ in $E_M$, whence $v_1=\ldots=v_d=0$
by Claim 3. This verifies the uniqueness. Next, we have
$\de(M)=\de(M^HA)=(M^H\ot1)\de(A)\sbs E_M$. Given $u\in M$, there exist
therefore $v_1,\ldots,v_d\in M$ such that
$$
\de u=\textsum\,(v_i\ot1)\cdot\de(a_i).\eqno(4.4)
$$
Applying $\id\ot\ep$ to both sides of this equality, we get $u=\sum v_ia_i$ since
$(\id_X\ot\ep)\circ\de=\id_X$ for both $X=M$ and $X=A$. Applying $\de\ot\id_H$
and $\id_M\ot\De$ to both sides of (4.4), and taking into account the identity
$(\de\ot\id)\circ\de=(\id\ot\De)\circ\de$, we get
$$
\textsum\,(\de v_i\ot1)\cdot(\de\ot\id)\de a_i=
\textsum\,(v_i\ot1\ot1)\cdot(\de\otimes\id)\de a_i\eqno(4.5)
$$
in $M\ot H\ot H$. Here we consider $M\ot H\ot H$ as a right $A\ot H\ot
H$-module in a natural way. Let $W$ denote the $A$-module $M\ot H$ on which $A$
operates via the ring homomorphism $\de:A\to A\ot H$. Then $A\ot H$ operates
on $W\ot H\cong M\ot H\ot H$ via the ring homomorphism $\de\ot\id_H$, and (4.5)
can be rewritten as an equality $\sum(w_i\ot1)\cdot\de(a_i)=0$ in $E_W$ where
$w_i=\de v_i-v_i\ot1\in W$. Claim 3 yields $\de v_i=v_i\ot1$, that is, $v_i\in
M^H$ for all $i$. Thus $u$ has the required form, and Claim 4 is proved.

Note that $A$ can be regarded as an object of $\M_A^H$. Moreover, $A\in\M'$
since $1\in A^H$. We see, in particular, that each element $a\in A$ can be uniquely
written as $\sum c_ia_i$ with $c_1,\ldots,c_d\in A^H$. In other words, $A$ is
a free $A^H$-module with $a_1,\ldots,a_d$ as its basis.

Now $M^H\otimes_{A^H}A\cong M^H\oplus\ldots\oplus M^H$ ($d$ copies), and the
restriction of $\Psi_M$ to the $i$th summand is given by the map $v\mapsto
va_i$. Claim 4 shows that $\Psi_M$ is bijective.

Suppose that $N$ is any $A^H$-module and $M=N\otimes_{A^H}A$. The composite
$$
N\otimes_{A^H}A\lmapr6{\Phi_N\ot\id}M^H\otimes_{A^H}A\lmapr4{\Psi_M}M
$$
is then the identity transformation of $M$. Since $\Psi_M$ is bijective, so
too is $\Phi_N\ot\id$. Since $A$ is a free $A^H$-module, $\Phi_N$ is bijective
as well. The proof of (i) and (ii) is now complete.

If $I$ is an $H$-costable ideal of $A$, then $A/I\in\M'$. As the canonical map
$A\to A/I$ is an epimorphism in $\M'$, the corresponding map $A^H\to(A/I)^H$ is an
epimorphism of $A^H$-modules by (ii). Thus $(A/I)^H\cong A^H/I^H$.
In the commutative diagram
$$
\diagram{
0&\mapr{}&I^H\otimes_{A^H}A&\mapr{}&A^H\otimes_{A^H}A&\mapr{}&A^H/I^H\otimes_{A^H}A&\mapr{}&0\cr
\diagramskip
&&&&\mapd{\Psi_A}{}&&\mapd{\Psi_{A/I}}{}&&\cr
\diagramskip
0&\lmapr4{}\hidewidth&I&\hidewidth\lmapr8{}\hidewidth&A&\hidewidth\lmapr8{}\hidewidth&A/I&\hidewidth\lmapr5{}&0\cr
}
$$
the rows are exact by freeness of $A$ over $A^H$ and the vertical arrows are
bijective by (ii). Hence $\Psi_A$ induces a bijection $I^H\otimes_{A^H}A\cong
I$, which shows that $I=I^HA$. Thus $I$ is an $\M'$-subobject of $A$. By (ii)
the $\M'$-subobjects of $A$ are in a bijective correspondence with the ideals
of $A^H$, and (iii) follows.
\endproof

\proclaim
Corollary 4.4.
If $A$ is $H$-reduced and $\p\in\Spec A$ is $H$-regular then there is an
isomorphism of $k(\p)$-algebras $O(\p)\cong k(\p)\ot_{k(\q)}A_\q/\q A_\q$
where $\q=\p\cap A^H$.
\endproclaim

\Proof.
In view of Lemma 3.2 we may replace $A$ with $A_\q$ and $\p$ with $\p A_\q$.
So we may assume that $A^H$ is local and $\q$ is its maximal ideal.
Then all prime ideals of $A$ are $H$-regular by Lemma 3.6, so that
the hypotheses of Theorem 4.3 are fulfilled. It follows that $A$ is a free
$A^H$-module of rank $d=\dim O(\p)$. Then $\dim_{k(\q)}A/\q A=d$. By Lemma 3.1
$\q\sbs\ker\de_\p$. Therefore $\de_\p$ extends to a homomorphism of
$k(\p)$-algebras $\ph:k(\p)\ot_{k(\q)}A/\q A\to O(\p)$ which is clearly
surjective. Comparing the dimensions, we deduce that $\ph$ is bijective.
\endproof

Consider the map $\ga:A\ot A\to A\ot H$ such that
$a\ot b\mapsto(a\ot1)\cdot\de(b)$. One says that $A$ is an {\it$H$-Galois
extension of $A^H$} if $\ga$ is surjective \cite{Chase}, \cite{Kr81}. In this
case $\ga$ induces a bijection $A\ot_{A^H}A\cong A\ot H$.

\proclaim
Proposition 4.5.
$A$ is an $H$-Galois extension of $A^H$ if and only if $O(\p)=k(\p)\ot H$ for
all $\p\in\Spec A$, if and only if $O(\p)=k(\p)\ot H$ for all maximal ideals
$\p$ of $A$.
\endproclaim

\Proof.
The map $\ga$ is $A$-linear with respect to the actions of $A$ by
multiplications on the first tensorands. If $\ga$ is surjective then so too is
$k(\p)\ot_A\ga$ for every $\p\in\Spec A$. The latter can be identified with the map
$\ga_\p:k(\p)\ot A\to k(\p)\ot H$ such that $c\ot b\mapsto(c\ot1)\de_\p(b)$.
The image of $\ga_\p$ coincides with $O(\p)$, and therefore $\ga_\p$
is surjective if and only if $O(\p)=k(\p)\ot H$. Note that the $A$-module $A\ot H$
is finitely generated since so is the $K$-module $H$. By \cite{Bou, Ch.~II,
\S3, Prop. 11} $\ga$ is surjective whenever the maps $\ga_\p$ are surjective for
all maximal ideals of $A$.
\endproof

\Remark.
If $A$ is $H$-Galois then $\dim O(\p)=\rk_{\p\cap K}H$. The function
$\p\mapsto\dim O(\p)$ is therefore locally constant in this case. In fact the
conclusions of Theorem 4.3 are valid for an $H$-Galois $A$ under much weaker
assumptions.  One needn't assume $A$ to be either $H$-reduced or even
commutative. If $A$ is $H$-Galois, then $A$ is finitely projective
as a left and as a right $A^H$-module \cite{Kr81, (1.7), (1.8)}. Moreover, the
functor $M\mapsto M^H$ establishes an equivalence between $\M_A^H$ and the
category of $A^H$-modules provided that $A$ is a faithfully flat left
$A^H$-module (this is automatic when $A$ is commutative) \cite{Sch90,
Th.~3.7}; similar results can be found in \cite{Doi89, (2.11)}, \cite{Ul82,
p.~661}. If $A$ and $H$ are commutative, then the condition that $A$ is
$H$-Galois means precisely that the finite group scheme $\Spec H$ operates
freely on $\Spec A$. The corresponding results on the quotient of this action
are contained in \cite{De, Ch.~III, \S12} and \cite{Mum70, Ch.~III, \S2}. Thus
Theorem 4.3 generalizes these classical results.
\endremark

\proclaim
Lemma 4.6.
Suppose that $A$ is $H$-reduced and Noetherian. Let $\p\in\Spec A$ be
$H$-regular, and let $a_1,\ldots,a_d\in A$ be elements such that
$\de_\p a_1,\ldots,\de_\p a_d$ form a basis for $O(\p)$ over $k(\p)$.
Then there exist $s\in A^H\setm\p$ and an integer $n>0$ such that
$s^nA\sbs A^Ha_1+\ldots+A^Ha_d$.
\endproclaim

\Proof.
Let $s$ be as in Lemma 3.4(iii). Put
\vadjust{\vskip-8pt}
$$
B=(A\ot1)\cdot\de(A)\qquad{\fam0and}\qquad
E=\textsum_{i=1}^d\,(A\ot1)\cdot\de(a_i)
$$
so that $E\sbs B\sbs A\ot H$. Both $B$ and $E$ are $A$-submodules of $A\ot H$
with respect to the action of $A$ by left multiplications on the first
tensorand. Localizing at $s$, we obtain a chain of $A_s$-modules $E_s\sbs
B_s\sbs A_s\ot H$. Clearly $B_s=(A_s\ot1)\cdot\de(A_s)$ and
$E_s=\sum_{i=1}^d(A_s\ot1)\cdot\de(a'_i)$ where $a'_i$ denotes the image of
$a_i$ in $A_s$. The $H$-comodule algebra $A_s$ fulfills the hypotheses under
which Claim 2 in the proof of Theorem 4.3 was verified. This claim shows that
$B_s$ is an $A_s$-module direct summand of $A_s\ot H$ and $B_s$ is a free
$A_s$-module with $\de(a'_1),\ldots,\de(a'_d)$ as its basis. In particular,
$B_s=E_s$. Now $B$ is a finitely generated $A$-module since so is $A\ot H$ and
$A$ is Noetherian. Hence there exists an integer $m>0$ such that
$(s^m\ot1)\cdot B\sbs E$.

Let $F$ be a free $A$-module with a basis $e_1,\ldots,e_d$. Consider the
$A$-module homomorphism $\ph:F\to A\ot H$ such that $e_i\mapsto\de(a_i)$ for
each $i=1,\ldots,d$, and put $R=\ker\ph$. The localization of $\ph$ at $s$
is an isomorphism of $F_s$ onto $B_s$, whence $R_s=0$. The $A$-module $R$ is
finitely generated as $A$ is Noetherian. It follows that there exists an
integer $l>0$ such that $s^lR=0$. Taking $n=m+l$, we will show that the
conclusion of the lemma is fulfilled.

Put $T=A\ot H$ considered as a ring extension of $A$ by means of the ring
homomorphism $\de:A\to A\ot H$. We claim that $T$ is a projective right
$A$-module with respect to $\de$. The $K$-linear transformation $\nu$ of
$A\ot H$ such that $a\ot h\mapsto(1\ot h)\cdot\de a$ for $a\in A$ and $h\in
H$ is invertible (the assignment
$a\ot h\mapsto(1\ot h)\cdot(\id\ot\si^{-1})(\de a)$ defines the inverse
transformation). As $\nu(ab\ot h)=\nu(a\ot h)\cdot\de b$ for all $a,b\in A$
and $h\in H$, one sees that $\nu$ is an isomorphism between the two $A$-module
structures on $A\ot H$ obtained via ring homomorphisms $b\mapsto b\ot1$ and
$\de$, respectively. As $H$ is $K$-projective, $A\ot H$ is projective with
respect to the first of these two module structures, hence also with respect
to the second one.

By projectivity of $T$ over $A$ we obtain an exact sequence of $T$-modules
$$
0\to T\ot_AR\to T\ot_AF\mapr\psi T\ot_A(A\ot H)\cong T\ot H\cong A\ot H\ot H
$$
where $\psi=\id\ot\ph$ and the isomorphisms shown are such that
$y\ot z$ in $T\ot_A(A\ot H)$ goes to $(y\ot1)\cdot(\de\ot\id)(z)$ in
$A\ot H\ot H$ where $y,z\in A\ot H$. As $\de s=s\ot1$ is in the center of $T$,
we get
$$
(s^l\ot1)T\ot_AR=T\de(s^l)\ot_AR=T\ot_As^lR=0.
$$
Note that $\psi(1\ot e_i)=(\de\ot\id)(\de a_i)$ in $A\ot H\ot H$ for each
$i=1,\ldots,d$. Given $x_1,\ldots,x_d\in T$, we deduce that
$$
\textsum\,(x_i\ot1)\cdot(\de\ot\id)(\de a_i)=0\Rar\textsum x_i\ot
e_i\in\ker\psi\Rar(s^l\ot1)x_i=0.\eqno(4.6)
$$

Let now $a\in A$ be any element, and put $b=s^ma$. Then
$\de b=(s^m\ot1)\cdot\de a\in E$. We can write therefore
$\de b=\sum\,(c_i\ot 1)\cdot\de a_i$ for some $c_1,\ldots,c_d\in A$. Applying
$\id\ot\ep$ to both sides of this equality, we get $b=\sum c_ia_i$.
Applying $\de\ot\id_H$ and $\id_A\ot\De$, we get
$$
\textsum\,(\de c_i\ot1)\cdot(\de\ot\id)\de a_i=
\textsum\,(c_i\ot1\ot1)\cdot(\de\otimes\id)\de a_i
$$
in $A\ot H\ot H$ (cf.\ the proof of Claim 4 in Theorem 4.3). Now (4.6) shows
that $(s^l\ot1)\cdot(\de c_i-c_i\ot1)=0$ for each $i=1,\ldots,d$. The last
equalities can be rewritten as $\de(s^lc_i)=s^lc_i\ot1$, which shows that
$s^lc_i\in A^H$ for each $i$. We deduce that $s^na=s^lb=\sum\,s^lc_ia_i$ is of
required form.
\endproof

\section
5. Semisimple stabilizer subalgebras and total integrals

For each $\p\in\Spec A$ consider the left coideal subalgebra $\St(\p)\sbs k(\p)\ot H^*$
corresponding to the right coideal subalgebra $O(\p)\sbs k(\p)\ot H$ (see Proposition 1.6).
We call $\St(\p)$ the {\it stabilizer subalgebra} associated with $\p$. This name is
justified by the next lemma. We may regard $k(\p)\ot A$ as an $H$-comodule algebra with
respect to the map $\id\ot\de$ or as a comodule algebra for the Hopf algebra
$k(\p)\ot H$ over $k(\p)$. Then $k(\p)\ot A$ is also a module over $k(\p)\ot H^*$.

\proclaim
Lemma 5.1.
Put $A'=k(\p)\ot A$, and let $\al_\p':A'\to k(\p)$ be the homomorphism of
$k(\p)$-algebras which extends the canonical ring homomorphism $\al_\p:A\to
k(\p)$. Then $\p'=\ker\al_\p'$ is a maximal ideal of $A'$ such that
$k(\p')\cong k(\p)$ and $\St(\p')\cong\St(\p)$. Furthermore, $\St(\p)$ is the
largest left coideal of $k(\p)\ot H^*$ which leaves $\p'$ stable.
\endproclaim

\Proof.
Clearly $A'/\p'\cong k(\p)$ so that $\p'$ is a maximal ideal of $A'$ with residue field
$k(\p)$. The map $\ga:A\to A'$ defined by the rule $a\mapsto1\ot a$ for
$a\in A$ is a homomorphism of $H$-comodule algebras. Since $\al'_\p\circ\ga=\al_\p$, we
have $\ga^{-1}(\p')=\ker\al_\p=\p$ and $\ga$ induces an isomorphism
$k(\p)\cong k(\p')$. Since $k(\p)\ot1\sbs A'^H$, Lemma 3.2 yields an isomorphism
$O(\p)\cong O(\p')$. Then $\St(\p)\cong\St(\p')$ as well. To complete the proof
of the lemma we may pass to base ring $k(\p)$ replacing $\p$, $A$, $H$ with
$\p'$, $A'$, $k(\p)\ot H$, respectively. Thus it suffices to consider the case where
$K$ is a field and $\p$ is a maximal ideal of $A$ such that $k(\p)\cong K$.

In this case $O(\p)=\de_\p(A)$. Since $\ep\circ\de_\p=\al_\p$ as shown in
diagram (3.1), we see that $\de_\p(a)\in O(\p)^+$ for $a\in A$ if and only if
$a\in\ker\al_\p=\p$. Hence $O(\p)^+=\de_\p(\p)$. For $\p$ to be stable under
the action of $\xi\in H^*$ it is necessary and sufficient that
$(\al_\p\circ L_\xi)(\p)=0$ where $L_\xi$ is the transformation of $A$ shown
in (1.2). Since
$$
\al_\p\circ L_\xi=(\al_\p\ot\xi)\circ\de=\xi\circ\de_\p,
$$
this condition can be rewritten as $\xi\bigl(O(\p)^+\bigr)=0$. Suppose that $C\sbs H^*$ is a
left coideal consisting of linear functions which vanish on $O(\p)^+$. If
$\xi\in C$ and $h\in H$ then $\xi\leftharpoonup h\in C$ where
$\xi\leftharpoonup h\in H^*$ is defined by the rule $(\xi\leftharpoonup
h)(g)=\xi(hg)$ for $g\in H$. Hence all $\xi\in C$ vanish on $H\cdot O(\p)^+$ and so
belong to $\St(\p)$.
\endproof

\Remark.
Note that in the situation of Proposition 1.6
$\dim C=\dim H/HB^+=\rk_BH=\dim H/\dim B$. Applying this fact to the orbital
subalgebras, we deduce that
$$
\dim O(\p)\cdot\dim\St(\p)=\dim k(\p)\ot H=\rk_{\p\cap K}H.
$$
\endremark

One may regard $H$ as a right $H$-comodule via $\De$. Any $H$-comodule
homomorphism $\ph:H\to A$ is called an {\it integral}. If, in addition,
$\ph(1)=1$ then $\ph$ is a {\it total integral}. As was observed by Doi \cite{Doi85}
the comodule algebras admitting a total integral enjoy remarkable properties.
An $H$-comodule $W$ is called {\it relatively injective} if, whenever $U$ is
an $H$-costable $K$-module direct summand of an $H$-comodule $V$, every
$H$-comodule map $U\to W$ can be extended to an $H$-comodule map $V\to W$.

\proclaim
Proposition 5.2.
The following conditions are equivalent{\rm:}

\item(i)
All objects $M\in\M_A^H$ are relatively injective $H$-comodules{\rm;}

\item(ii)
$A$ is a relatively injective $H$-comodule{\rm;}

\item(iii)
There exists a total integral $\ph:H\to A${\rm;}

\item(iv)
There exist natural $K$-linear retractions $\tr_M:M\to M^H$, defined for each
$M\in\M_A^H${\rm;}

\item(v)
The functor $M\mapsto M^H$ is exact on $\M_A^H${\rm;}

\item(vi)
$A$ is a projective object of $\M_A^H$.

\endproclaim

\Proof.
Equivalences (i)$\Lrar$(ii)$\Lrar$(iii) are proved in \cite{Doi85, (1.6)}.

(iii)$\Rar$(iv)
Denote by $\tr_M$ the composite
$$
M\lmapr2{\de}M\ot H\lmapr4{\id\ot\si}M\ot H\lmapr4{\id\ot\ph}M\otimes A\lmapr2{}M
$$
where the last map is afforded by the $A$-module structure on $M$. As was shown in
\cite{Doi90, \S1} $\tr_M(M)\sbs M^H$ and all required properties are fulfilled.

(iv)$\Rar$(v)
The functor $M\mapsto M^H$ is clearly left exact. Suppose that $\xi:M\to N$ is
an epimorphism in $\M_A^H$. Given $u\in N^H$, take $v\in M$ such that
$\xi(v)=u$. Then $\tr_M(v)\in M^H$ and
$\xi\bigl(\tr_M(v)\bigr)=\tr_N\bigl(\xi(v)\bigr)=u$. Thus $\xi$ induces a
surjection $M^H\to N^H$.

(v)$\Lrar$(vi)
This is clear since every morphism $A\to M$ in $\M_A^H$ is given by the rule $a\mapsto va$ where
$v\in M^H$.

(v)$\Rar$(iii)
If $V,W$ are two $K$-modules and $V$ is finitely projective,
then there is a canonical bijection $\al_{VW}:V^*\ot W\to\Hom_K(V,W)$ where
$V^*=\Hom_K(V,K)$. Suppose that $V,W$ are, moreover, $H$-comodules. Define an
$H$-comodule structure on $V^*$ such that $\de\eta\in V^*\ot H$, where
$\eta\in V^*$, corresponds under bijection $\al_{VH}$ to the composite
$K$-linear map
$$
V\lmapr2{\de}V\ot H\lmapr4{\eta\ot\si^{-1}}K\ot H\cong H
$$
(note that the antipode $\si$ of $H$ is bijective by \cite{Par71, Prop.~4}).
If $V^*\ot W$ is equipped with the tensor product of two comodule structures,
then $\al_{VW}$ maps $(V^*\ot W)^H$ bijectively onto the set $\Com(V,W)$ of
all $H$-comodule homomorphisms $V\to W$. As a special case we have
$V^*\ot A\in\M_A^H$ where $A$ operates by multiplications on the second
tensorand. Now $K$ is a $K$-module direct summand of $H$ as $\ep:H\to K$ is a
retraction. Moreover, $K$ is an $H$-subcomodule of $H$. Then the restriction
map $H^*\to K^*$ is an epimorphism of $H$-comodules which gives rise to an
epimorphism $H^*\ot A\to K^*\ot A$ in $\M_A^H$. Taking the invariants, we
deduce that the restriction map $\Com(H,A)\to\Com(K,A)$ is surjective. In
particular, the $H$-comodule map $K\to A$ such that $1\mapsto1$ extends
to an integral $\ph:H\to A$.
\endproof

\Remark.
The commutativity of $A$ is not needed in this result. The assumption that $H$
is a finitely projective $K$-module was used only in the proof of (v)$\Rar$(iii).
\endremark

\proclaim
Proposition 5.3.
If the algebra $\St(\p)$ is semisimple for some $\p\in\Spec A$ then there exists
an integral $\ph:H\to A$ such that $\ph(1)\notin\p$. If the algebras $\St(\p)$ are
semisimple for all maximal ideals $\p$ of $A$ then there exists a total integral
$\ph:H\to A$.
\endproclaim

\Proof.
(i) Consider first a {\it Preliminary Step} in which we assume that $K$ is a
field and $\p$ is a maximal ideal of $A$ such that $A=K+\p$. Then
$k(\p)\cong K$, and one has $O(\p)=\de_\p(A)$.
Let $I=\ker\de_\p$. The $H$-comodule algebras $O(\p)$ and $A/I$ are
isomorphic, whence the category $\M_{A/I}^H$ is equivalent to the category
of $\St(\p)$-modules by Proposition 1.6. As $\St(\p)$ is semisimple, all
$\St(\p)$-modules are projective. It follows that so are all objects of
$\M_{A/I}^H$ too. In particular, $A/I$ is projective in $\M_{A/I}^H$. Then
$A/I$ is an injective $H$-comodule by Proposition 5.2. Hence $A/I$ is an
injective $H^*$-module. Since $H^*$ is a Frobenius algebra, $A/I$ is also a
projective $H^*$-module. It follows that the canonical epimorphism of
$H^*$-modules $A\to A/I$ splits, and so there exists an $H$-subcomodule
$V\sbs A$ such that $A=V\oplus I$. Now $V$ is an injective $H$-comodule as
it is isomorphic to $A/I$. Denote by $v\in V$ the element which projects to
1 in $A/I$. Clearly $v\in V^H$. Note that $K$ is a trivial $H$-subcomodule
of $H$. By injectivity of $V$ the $H$-comodule homomorphism $K\to V$ such
that $1\mapsto v$ extends to an $H$-comodule homomorphism $\ph:H\to V$ which
can be regarded as an integral $H\to A$. We have $\ph(1)=v$. By Lemma 3.1
$I\sbs\p$. It follows that $v\notin\p$, and the required property of $\ph$ is fulfilled.

Consider now the {\it General Case}. Let $\p$ be given, and put $H'=k(\p)\ot H$.
The hypotheses of the Preliminary Step are fulfilled for the $H'$-comodule algebra
$A'=k(\p)\ot A$ over $k(\p)$ and its maximal ideal $\p'$ defined in Lemma 5.1.
It follows that there exists an integral $\psi:H'\to A'$ such that
$\psi(1)\notin\p'$. Taking the composite with the map $H\to H'$ defined by the
rule $h\mapsto1\ot h$, we obtain an integral $\chi:H\to A'$ satisfying
$\chi(1)\notin\p'$. Note that the $H$-comodule homomorphisms are precisely the $H^*$-module
homomorphisms. By Lemma 1.3 the $H^*$-module $H$ is finitely projective. Hence
$$
\Hom_{H^*}(H,A')\cong k(\p)\ot\Hom_{H^*}(H,A).
$$
If $\ph:H\to A$ is an integral, then $1\ot\ph$ corresponds under this isomorphism to the
integral $\ga\circ\ph:H\to A'$ where $\ga:A\to A'$ is defined by the rule $a\mapsto1\ot a$
for $a\in A$. There exists therefore $\ph$ such that $(\ga\circ\ph)(1)\notin\p'$.
Since $\ga^{-1}(\p')=\p$, we get $\ph(1)\notin\p$.

(ii) If $a\in A^H$ and $\ph,\psi:H\to A$ are two integrals then $a\ph$ and
$\ph+\psi$ are also integrals. In other words, the set $\Com(H,A)$ of all
integrals $H\to A$ is an $A^H$-module in a natural way. One sees that the map
$\Com(H,A)\to A$ such that $\ph\mapsto\ph(1)$ is a homomorphism of
$A^H$-modules. Its image $J$ is therefore an ideal of $A^H$. If $\St(\p)$ is
semisimple for some $\p\in\Spec A$ then $J\not\sbs\p$ by (i). It follows that
$J=A$ as long as the algebras $\St(\p)$ are semisimple for all maximal ideals
of $A$. In this case $1\in J$, and we are done.
\endproof

\Remark.
Conversely, if there exists an integral $\ph:H\to A$ such that $\ph(1)\notin\p$,
then it is possible to prove, using Proposition 1.6, that all left $k(\p)\ot
H^*$-modules are semisimple $\St(\p)$-modules. I don't know if this suffices
to assert that $\St(\p)$ is semisimple. Koppinen \cite{Kop93, section 5} shows
the semisimplicity of coideal subalgebras in finite dimensional Hopf algebras under
additional assumptions. If, however, $H$ is commutative, then $\St(\p)$ is a Hopf
subalgebra, and the conclusion is true. In this case one can add another equivalent
condition to the list in Proposition 5.2:

\item(vii)
{\it the algebras $\St(\p)$ are semisimple for all $\p\in\Spec A$.}
\endremark

We say that the Hopf algebra $H$ is {\it geometrically cosemisimple} if for
every ring homomorphism $K\to E$ into a field $E$ the Hopf algebra $E\ot H$
over $E$ is cosemisimple, that is, the dual Hopf algebra $E\ot H^*$ is semisimple.

\proclaim
Corollary 5.4.
$H$ is geometrically cosemisimple if and only if $H^*$ contains a right
integral $\ph$ such that $\ph(1)=1$.
\endproclaim

\Proof.
Apply Proposition 5.3 to $A=K$ equipped with the trivial coaction of $H$, so
that $K^H=K$. We have $O(\p)=k(\p)\ot1$, hence $\St(\p)=k(\p)\ot H^*$ for
each $\p\in\Spec A$. If $H$ is geometrically cosemisimple then the algebras
$k(\p)\ot H^*$ are semisimple. In this case there exists a total integral
$\ph:H\to K$ which is none other but a right integral in $H^*$. Conversely, if
$H^*$ has a right integral $\ph$ satisfying $\ph(1)=1$, then so too
does $E\ot H^*$ for any ring homomorphism $K\to E$ into a field $E$. By
\cite{Sw, Th.~5.1.8} $E\ot H^*$ is semisimple.
\endproof

As an example consider one special case: Suppose that $pK=0$ where $p$ is a
prime integer. Let $L$ be a $p$-Lie algebra over $K$ whose underlying
$K$-module is finitely projective. Denote by $u(L)$ the restricted universal
enveloping algebra of $L$. Its underlying $K$-module is finitely projective
as well \cite{De, Ch.~II, \S7, Cor.~3.7}. Put $H=u(L)^*$ so that $H^*\cong u(L)$. Every
homomorphism of $p$-Lie algebras $L\to\Der_KA$ into the $p$-Lie algebra of
$K$-linear derivations of $A$ gives rise to an action of $u(L)$ on $A$ which
makes $A$ into a $u(L)$-module algebra, hence into an $H$-comodule algebra. If
$\p\in\Spec A$ then $\St(\p)$ is a Hopf subalgebra of $k(\p)\ot u(L)\cong
u\bigl(k(\p)\ot L\bigr)$ since the latter Hopf algebra is cocommutative. Hence
$\St(\p)=u(L_\p)$ where $L_\p$ is a $[p]$-closed subalgebra of the finite
dimensional $p$-Lie algebra $k(\p)\ot L$ over $k(\p)$. In view of Lemma 5.1
$L_\p$ coincides with the stabilizer of $\p'$ with respect to the natural
action of $k(\p)\ot L$ on $k(\p)\ot A$.

\proclaim
Corollary 5.5.
In order that there exist a total integral $u(L)^*\to A$ it is necessary and
sufficient that $L_\p$ be a torus for every maximal ideal $\p$ of $A$.
\endproclaim

\Proof.
As was shown by Hochschild \cite{Hoch54} $u(L_\p)$ is semisimple if and only if
$L_\p$ is a torus. Hence Proposition 5.3 and the Remark following it applies.
\endproof

\section
6. Weakly reductive coactions

We say that $A$ is {\it weakly reductive} with respect to the coaction of
$H$ if for every $H$-costable ideal $I$ of $A$ the canonical map
$\psi:A^H\to(A/I)^H$ is surjective.

\proclaim
Proposition 6.1.
$A$ is {\it weakly reductive} with respect to coaction of $H$ in any of
the following three cases:

\item(a)
$A$ has invariant characteristic polynomials and $\rk_\q H$ is invertible in $A$ for every $\q\in\Spec K$,

\item(b)
there exists a total integral $H\to A$,

\item(c)
$A$ is $H$-reduced and $\St(\p)$ is semisimple for every maximal ideal $\p$
of $A$ which is not $H$-regular.

\endproclaim

\Proof.
Let $I$ be an $H$-costable ideal of $A$.

(a) Denote by $m$ the least common multiple of all $\rk_\q H$ with
$\q\in\Spec K$. Proposition 2.6 shows that $m\cdot(A/I)^H\sbs\im\psi$. The
hypothesis of (a) imply that $m^{-1}\in A$, hence $m^{-1}\in A^H$. It
follows that $\psi$ is surjective.

(b) The surjectivity of $\psi$ is a consequence of condition (v) in
Proposition 5.2.

(c) For $\psi$ to be surjective it is necessary and sufficient that so be its localizations
$\psi_\q$ at all maximal ideals $\q$ of $A^H$. By Lemma 1.1 we may identify
$\psi_\q$ with the canonical map $(A_\q)^H\to(A_\q/IA_\q)^H$. Replacing $A$
with $A_\q$, we may thus assume that the ring $A^H$ is local with a maximal
ideal $\q$. If there exists an $H$-regular maximal ideal of $A$ then all prime
ideals of $A$ are $H$-regular by Lemma 3.6, and the surjectivity of $\psi$ follows from
Theorem 4.3(ii). Otherwise the algebras $\St(\p)$ are semisimple for all maximal ideals
of $A$, and there exists a total integral $H\to A$ by Proposition 5.3. In this
case (b) applies.
\endproof

\Remark.
Suppose that $K$ is field. In this case the second part of (a) means that
$\chr K$ does not divide $\dim_KH$. Condition (b) is fulfilled, for instance,
when $H$ is cosemisimple.
\endremark

\proclaim
Theorem 6.2.
If $A$ is weakly reductive with respect to coaction of $H$ then{\rm:}

\item(i)
$A$ is integral over $A^H$.

\item(ii)
If $\p\in\Spec A$ and $\q=\p\cap A^H$ then $k(\p)$ is a finite extension of
$k(\q)$ of degree $[k(\p):k(\q)]\le\dim O(\p)$.

\item(iii)
If $A$ is Noetherian then so too is $A^H$ and $A$ is a finite $A^H$-module.

\item(iv)
If $A$ is Cohen-Macaulay then so too is $A^H$.

\endproclaim

\Proof.
(i) Denote by $R$ the $H$-radical of $A$. As $A/R$ is $H$-reduced, it has invariant
characteristic polynomials by Theorem 2.5, and so $A/R$ is integral over $(A/R)^H$
by Proposition 2.7. Since $A$ is weakly reductive, we have $(A/R)^H=B/R$ where
$B=A^H+R$. It follows that $A$ is integral over $B$. All elements of $R$ are
nilpotent, hence integral over $A^H$. Then $B$ is integral over $A^H$, and so too is $A$.

(ii) Let $I$ be the largest $H$-costable ideal of $A$ such that $I\sbs\p$.
The $H$-comodule algebra $A'=A_\q/IA_\q$ is $H$-simple by Lemma 4.2.
Letting $\p'=\p A_\q/IA_\q$ and $\q'=\p'\cap A'^H$, we have $k(\p')\cong k(\p)$
and $k(\q')\cong k(\q)$. The second isomorphism is a consequence of the
surjectivity of the map $A^H\to(A/I)^H$. We may therefore replace $A$ with
$A'$ and so assume that $A$ is $H$-simple. Then the desired conclusion
follows from Lemma 4.1. Indeed, $A^H$ is a field, so that $\q=(0)$ and
$k(\q)\cong A^H$. We have also $\dim_{A^H}A=\dim O(\p)$ and $k(\p)\cong
A/\p$.

(iii) Suppose that $A$ is not a finitely generated $A^H$-module. Since
$A$ is Noetherian, $A$ contains an ideal $I$, maximal with respect to the
properties that $I$ is $H$-costable and $A/I$ is not a finitely generated
$(A/I)^H$-module. Replacing $A$ with $A/I$, we may assume that $A/J$ is a
finitely generated $(A/J)^H$-module for every $H$-costable nonzero ideal $J$
of $A$. Since the canonical map $A^H\to(A/J)^H$ is surjective, for any such
a $J$ there exists a finitely generated $A^H$-submodule $F\sbs A$ which is
mapped onto the whole $A/J$ under the projection $A\to A/J$, so that $A=F+J$.

Denote by $R$ the $H$-radical of $A$. Suppose that $R\ne0$. Then there
exists a finitely generated $A^H$-submodule $F\sbs A$ such that $A=F+R$. Let
$X$ be a finite set of generators for the $A^H$-module $F$. Since $A$ is
Noetherian, the ideal $R$ of $A$ has a finite set of generators, say $Y$.
Denote by $B$ the $A^H$-subalgebra of $A$ generated by $X\cup Y$. As $F\sbs B$ by
construction, we have $A=B+R$. Since $Y\sbs B$, we have also $R=(B+R)Y\sbs
B+R^2$. An easy induction shows that $R^{m-1}=(R^{m-1}\cap B)+R^m$ and
$A=B+R^m$ for all $m>0$. As $R$ is a finitely generated nil ideal of $A$,
there exists $m$ such that $R^m=0$. Thus $A=B$, and so $A$ is finitely generated
as an algebra over $A^H$. By (i) $A$ is also a finitely generated $A^H$-module
\cite{Bou, Ch.~V, \S1, Prop.~4}.

Suppose now that $R=0$, that is, $A$ is $H$-reduced. By Proposition 3.5 $A$
contains an $H$-regular prime ideal $\p$. Let $a_1,\ldots,a_d\in A$, $\,s\in
A^H$ and $n>0$ be as in Lemma 4.6. The ideal $s^nA$ of $A$ is $H$-costable
and nonzero. Hence there exists a finitely generated $A^H$-submodule $F\sbs A$
such that $A=F+s^nA$. From Lemma 4.6 it follows that $A=F+\sum A^Ha_i$.
Thus in both cases we have arrived at a contradiction.

We can conclude that $A$ is a finitely generated $A^H$-module. By Eakin's
theorem \cite{Ea}, \cite{Ma, Appendix}, \cite{Na} $A^H$ is Noetherian.

(iv) By (iii) $A^H$ is Noetherian. We have to prove that $A^H_\q$ is
Cohen-Macaulay for every $\q\in\Spec A^H$. Replacing $A$ with $A_\q$, we may
thus assume that the ring $A^H$ is local with a maximal ideal $\q$. If
$\Spec A^H=\{\q\}$ then $A^H$ is Artinian, and we are done. Otherwise, let
$\p_1,\ldots,\p_n$ be all minimal prime ideals of $A$. By the going-down
property in Theorem 3.3 each $\p_i\cap A^H$ is a minimal prime ideal of
$A^H$, whence $\q\not\sbs\p_i$. Then $\q\not\sbs\bigcup\p_i$ \cite{Bou, Ch.~II,
\S1, Prop.~2}. The set of zero divisors in $A$ coincides with $\bigcup\p_i$
\cite{Ma, (16.C)}. It follows that $\q$ contains an element $x$ which is not a zero
divisor in $A$. The factor ring $A/xA$ is Cohen-Macaulay of smaller Krull dimension
than $A$. Furthermore, $(A/xA)^H\cong A^H/(xA\cap A^H)$ since $A$ is weakly reductive.
Suppose that $a\in A$ is an element such that $xa\in A^H$. Then $(x\ot1)\de
a=\de(xa)=xa\ot1$. Since $H$ is a flat $K$-module, $x\ot1$ is not a zero
divisor in $A\ot H$. So it follows that $\de a=a\ot1$, i.e., $a\in A^H$.
Thus $(A/xA)^H\cong A^H/xA^H$. If the ring $A^H/xA^H$ is Cohen-Macaulay
then so too is $A^H$ since $x$ is not a zero divisor in $A^H$ \cite{Ma, (16.A)}.
So we may proceed by induction on $\dim A$.
\endproof

\Remarks.
If there exists a total integral $H\to A$ then $A^H$ is an $A^H$-module direct
summand of $A$. In this case (iv) follows from the Hochster-Eagon theorem
\cite{Br, Th.~6.4.5}, \cite{Ho71, Prop.~12}.

One may also note that all results of section 3 are valid for a weakly
reductive $A$ without the assumption that $A$ has invariant characteristic
polynomials. Indeed, $\Spec A$ and $\Spec A^H$ are homeomorphic to $\Spec A/R$
and $\Spec(A/R)^H$, respectively, where $R$ stands for the $H$-radical of $A$.
So the proofs are obtained by passing to the $H$-reduced algebra $A/R$.
\endremark

Suppose that $pK=0$ where $p$ is a prime integer and $L$ is a finitely $K$-projective
$p$-Lie algebra over $K$ operating on $A$ via $K$-linear derivations. Taking into
account Corollary 5.5, we get

\proclaim
Corollary 6.3.
If $A$ is Cohen-Macaulay and $L_\p$ is a torus for every maximal ideal $\p$ of
$A$ then the subring of $L$-invariants $A^L\sbs A$ is Cohen-Macaulay.
\endproclaim

This generalizes \cite{Ar86, (7.1)} and \cite{Ar87, (3.1)} where the
invariants of a single derivation were considered under more restrictive
assumptions on $A$.
\references
\nextref
Ar86
\auth
A.G.,Aramova;L.L.,Avramov;
\endauth
\paper{Singularities of quotients by vector fields in characteristic $p$}
\journal{Math. Ann.}
\Vol{273}
\Year{1986}
\Pages{629--645}

\nextref
Ar87
\auth
A.G.,Aramova;
\endauth
\paper{Reductive derivations of local rings of characteristic $p$}
\journal{J. Algebra}
\Vol{109}
\Year{1987}
\Pages{394--414}

\nextref
BouA
\auth
N.,Bourbaki;
\endauth
\book{Alg\`ebre, Chapitres 1--3}
\publisher{Hermann}
\Year{1970};
\book{Chapitres 4--7}
\publisher{Masson}
\Year{1981}.

\nextref
Bou
\auth
N.,Bourbaki;
\endauth
\book{Commutative Algebra}
\publisher{Springer}
\Year{1989}.

\nextref
Br
\auth
W.,Bruns;J.,Herzog;
\endauth
\book{Cohen--Macaulay Rings}
\publisher{Cambridge University Press}
\Year{1993}.

\nextref
Chase
\auth
S.U.,Chase;M.,Sweedler;
\endauth
\book{Hopf Algebras and Galois Theory}
\bookseries{Lecture Notes in Mathematics}
\Vol{97}
\publisher{Springer}
\Year{1969}.

\nextref
De
\auth
M.,Demazure;P.,Gabriel;
\endauth
\book{Groupes Alg\'ebriques I}
\publisher{Masson}
\Year{1970}.

\nextref
Doi85
\auth
Y.,Doi;
\endauth
\paper{Algebras with total integrals}
\journal{Comm. Algebra}
\Vol{13}
\Year{1985}
\Pages{2137--2159}

\nextref
Doi89
\auth
Y.,Doi;M.,Takeuchi;
\endauth
\paper{Hopf--Galois extensions of algebras, the Miyashita--Ulbrich action, and Azumaya algebras}
\journal{J. Algebra}
\Vol{121}
\Year{1989}
\Pages{488--516}

\nextref
Doi90
\auth
Y.,Doi;
\endauth
\paper{Hopf extensions of algebras and Maschke type theorems}
\journal{Isr. J.~Math.}
\Vol{72}
\Year{1990}
\Pages{99--108}

\nextref
Ea
\auth
P.M.,Eakin Jr.;
\endauth
\paper{The converse to a well known theorem on Noetherian rings}
\journal{Math. Ann.}
\Vol{177}
\Year{1968}
\Pages{278--282}

\nextref
Hab75
\auth
W.J.,Haboush;
\endauth
\paper{Reductive groups are geometrically reductive}
\journal{Ann. Math.}
\Vol{102}
\Year{1975}
\Pages{67--83}

\nextref
Hoch54
\auth
G.P.,Hochschild;
\endauth
\paper{Representations of restricted Lie algebras of characteristic $p$}
\journal{Proc. Amer. Math. Soc.}
\Vol{5}
\Year{1954}
\Pages{603--605}

\nextref
Ho71
\auth
M.,Hochster;J.A.,Eagon;
\endauth
\paper{Cohen--Macaulay rings, invariant theory, and the generic perfection of determinantal loci}
\journal{Amer. J. Math.}
\Vol{93}
\Year{1971}
\Pages{1020--1058}

\nextref
Kop93
\auth
M.,Koppinen;
\endauth
\paper{Coideal subalgebras in Hopf algebras: Freeness, integrals, smash products}
\journal{Comm. Algebra}
\Vol{21}
\Year{1993}
\Pages{427--444}

\nextref
Kr81
\auth
H.F.,Kreimer;M.,Takeuchi;
\endauth
\paper{Hopf algebras and Galois extensions of an algebra}
\journal{Indiana Univ. Math.~J.}
\Vol{30}
\Year{1981}
\Pages{675--692}

\nextref
Lar95
\auth
R.G.,Larson;D.E.,Radford;
\endauth
\paper{Semisimple Hopf algebras}
\journal{J. Algebra}
\Vol{171}
\Year{1995}
\Pages{5--35}

\nextref
Mas92
\auth
A.,Masuoka;
\endauth
\paper{Freeness of Hopf algebras over coideal subalgebras}
\journal{Comm. Algebra}
\Vol{20}
\Year{1992}
\Pages{1353--1373}

\nextref
Ma
\auth
H.,Matsumura;
\endauth
\book{Commutative Algebra, Second Edition}
\publisher{Benjamin}
\Year{1980}.

\nextref
Mo
\auth
S.,Montgomery;
\endauth
\book{Hopf algebras and Their Actions on Rings}
\bookseries{CBMS Regional Conference Series in Mathematics}
\Vol{82}
\publisher{American Mathematical Society}
\Year{1993}.

\nextref
Mum65
\auth
D.,Mumford;
\endauth
\book{Geometric Invariant Theory}
\bookseries{Ergebnisse der Mathematik und ihrer Grenzgebiete}
\Vol{34}
\publisher{Springer}
\Year{1965}.

\nextref
Mum70
\auth
D.,Mumford;
\endauth
\book{Abelian Varieties}
\publisher{Oxford University Press}
\Year{1970}.

\nextref
Na
\auth
M.,Nagata;
\endauth
\paper{A type of subrings of a noetherian ring}
\journal{J. Math. Kyoto Univ.}
\Vol{8}
\Year{1968}
\Pages{465--467}

\nextref
Par71
\auth
B.,Pareigis;
\endauth
\paper{When Hopf algebras are Frobenius algebras}
\journal{J. Algebra}
\Vol{18}
\Year{1971}
\Pages{588--596}

\nextref
Sch90
\auth
H.--J.,Schneider;
\endauth
\paper{Principal homogeneous spaces for arbitrary Hopf algebras}
\journal{Isr. J.~Math.}
\Vol{72}
\Year{1990}
\Pages{167--195}

\nextref
Sk02
\auth
S.,Skryabin;
\endauth
\paper{Invariants of finite group schemes}
\journal{J.~London Math. Soc.}
\Vol{65}
\Year{2002}
\Pages{339--360}

\nextref
Sw
\auth
M.E.,Sweedler;
\endauth
\book{Hopf Algebras}
\publisher{Benjamin}
\Year{1969}.

\nextref
Ul82
\auth
K.--H.,Ulbrich;
\endauth
\paper{Galoiserweiterungen von nicht--kommutativen Ringen}
\journal{Comm. Algebra}
\Vol{10}
\Year{1982}
\Pages{655--672}

\nextref
Wat
\auth
W.C.,Waterhouse;
\endauth
\paper{Geometrically reductive affine group schemes}
\journal{Arch. Math.}
\Vol{62}
\Year{1994}
\Pages{306--307}

\nextref
Zhu96
\auth
S.,Zhu;
\endauth
\paper{Integrality of module algebras over its invariants}
\journal{J. Algebra}
\Vol{180}
\Year{1996}
\Pages{187--205}

\endreferences
\bye